\documentclass[10pt,twoside]{article}
\usepackage{amsmath,amssymb,amsthm}
\usepackage{epsfig,graphicx,verbatim}
\usepackage{color}
\usepackage[bookmarks, bookmarksopen=true, bookmarksnumbered=true]{hyperref}
\numberwithin{equation}{section}
\newtheorem{thm}[equation]{Theorem}
\newtheorem{prop}[equation]{Proposition}

\newtheorem{lem}[equation]{Lemma}

\newcounter{mycount}
\newenvironment{romlist}{\begin{list}{\rm(\roman{mycount})}%
   {\usecounter{mycount}\labelwidth=1cm\itemsep 0pt}}{\end{list}}

\newenvironment{letlist}{\begin{list}{(\alph{mycount})}%
   {\usecounter{mycount}\labelwidth=1cm\itemsep 0pt}}{\end{list}}
\def\Wapq{{\mathcal W}_{a,p,q}}
\def\Rapq{{\mathcal R}_{a,p,q}}
\def\coWapq{\overline{\mathrm{co}\,{\mathcal W}_{a,p,q}}}
\def\Vapq{{\mathcal V}_{a,p,q}}
\def\coVapq{\overline{\mathrm{co}\,{\mathcal V}_{a,p,q}}}

\def\Ld{\mathbb L^d}
\def\Ed{{\mathbb E}^d}
\def\Zd{{\mathbb Z}^d}

\def\L2{{{\mathbb L}^2}}
\def\sF{\mathcal F}
\def\sG{\mathcal G}

\def\sS{\mathcal S}
\def\sH{\mathcal H}

\def\laa{{\lambda_1}}
\def\lab{{\lambda_2}}
\def\lai{{\lambda_i}}
\def\po{\psi,\omega}
\def\xy{{\langle x,y\rangle}}
\def\sT{{\mathcal T}}
\def\ZBC{Z^{\text{BC}}}
\def\ZRC{Z^{\text{RC}}}
\def\ZDRC{Z^{\text{DRC}}}
\def\ZBCP{Z^{\text{BCP}}}
\def\HH{{\mathcal H}}
\def\s{\sigma}
\def\ZI{Z^{\text{I}}}
\def\qq{\qquad}
\def\q{\quad}

\def\b{\beta}
\def\d{\delta}
\def\l{\lambda}
\def\t{\theta}
\def\rc{random-cluster}
\def\BC{Blume--Capel}
\def\BCP{Blume--Capel--Potts}
\def\ZZ{{\mathbb Z}}
\def\RR{{\mathbb R}}
\def\LL{{\mathbb L}}

\def\la{\langle}
\def\ra{\rangle}
\def\Si{\Sigma}
\def\Om{\Omega}
\def\om{\omega}
\def\De{\Delta}
\def\Th{\Theta}
\def\La{\Lambda}
\def\oo{\infty}
\def\drc{diluted-\rc}
\def\deg{\text{deg}}
\def\lest{\le_{\mathrm {st}}}

\def\Kc{K_{\mathrm c}}
\def\gest{\ge_{\mathrm {st}}}
\def\Var{\mathrm{var}}
\def\be{\begin{equation}}
\def\ee{\end{equation}}
\def\sm{\setminus}
\def\resp{respectively}
\def\mupq{\mu_{p,q}}
\def\EM{\Upsilon}

\def\pcsite{p_{\rm{c}}^{\rm{site}}}
\def\pcb{p_{\rm{c}}^{\rm{bond}}}
\def\Jc{J_{\rm{c}}}
\def\pic{\pi_{\rm{c}}}

\def\sB{\mathcal {B}}
\def\bzero{{\mathbf{0}}}
\def\bone{{\mathbf{1}}}

\def\pd{\partial}
\def\lra{\leftrightarrow}
\def\papq{\phi^\bone_{a,p,q}}
\def\Picvc{\Pi_{\text{\rm icvc}}}
\def\Piovc{\Pi_{\text{\rm iovc}}}
\def\Piec{\Pi_{\text{\rm iec}}}
\def\es{\varnothing}
\def\paptwo{\phi^\bone_{a,p,2}}

\def\ol#1{\overline{#1}}
\def\RC{{\rm RC}}
\def\P{{\rm P}}

\begin{document}
\title{Random-cluster representation of the \BC\ model}
\author{B.\ T.\ Graham, G.\ R.\ Grimmett\\
   {\small Statistical Laboratory, University of Cambridge,}\\
   {\small Wilberforce Road, Cambridge CB3 0WB, U.K.}}
\maketitle

\begin{abstract}
The so-called \drc\ model may be viewed
as a \rc\ representation
of the Blume--Capel model. 
It has three parameters, a vertex parameter $a$,
an edge parameter $p$, and a cluster weighting factor $q$.
Stochastic comparisons of measures are developed for the `vertex marginal'
when $q\in[1,2]$, and the `edge marginal'
when $q\in[1,\oo)$. Taken in conjunction with
arguments used earlier for the \rc\ model, these
permit a rigorous study of part of the phase diagram
of the Blume--Capel model.
\\

\noindent
{\bf Keywords} Blume--Capel model, Ising model, Potts model, \rc\ model,
first-order phase transition, tri-critical point.
\\

\noindent
{\bf Mathematics Subject Classification (2000)} 82B20, 60K35.
\end{abstract}

\section{Introduction}
The Ising model is one of the most studied models of statistical 
physics. It
has configuration space $\{-1,+1\}^V$ where $V$ is the vertex set
of the (finite) graph $G$ in question, and has Hamiltonian
\[
\HH(\s) = -J
\sum_{\xy} \s_x\s_y  
 - h \sum_{x\in V} \s_x, \qq \s\in\{-1,+1\}^V.
\]
The first summation is over all (unordered)
pairs of nearest neighbours, and $J\in[0,\oo)$, $h\in\RR$.
The Ising probability measure $\mu$ on $\{-1,+1\}^V$ is given by
\[
\mu(\s) = \frac 1 {\ZI} e^{-\b\HH(\s)},\qq \s\in\{-1,+1\}^V,
\]
where $\ZI$ is the appropriate normalizing constant.
Here, $\b=1/(kT)$ where $k$ is Boltzmann's constant and $T$ is
temperature.

It is standard that the Ising measure may be extended to a probability
measure on the configuration space associated with an infinite graph.
For physical and mathematical reasons, it is convenient that this
graph have a good deal of symmetry, and it is usual
to work with the $d$-dimensional hypercubic 
lattice $\ZZ^d$, where $d\ge 2$. In such a case, the model undergoes
a phase transition, and this is the main phenomenon
of interest in the theory. This transition is known
to be of second-order (continuous) when $d=2$ or $d \ge 4$,
and is believed to be of second-order when $d=3$ also.
See \cite{ABF,AF,FFS}.

The Ising model has two local states, namely $\pm 1$.
This may be generalized to any given number $q\in\{2,3,\dots\}$
of local states by considering the  
so-called Potts model introduced in 1952, see
\cite{POTTS}.  The Potts phase transition
is richer in structure than that
of the Ising model, in that it is of first-order (discontinuous) if $q$ is 
sufficiently large. See \cite{G-RC,Kot-S,LMMRS}.

In 1966, Blume introduced a variant of the Ising model,
see \cite{BLUME}, with the physical motivation of studying 
magnetization in Uranium Oxide, UO$_2$, at a
temperature of
about $30^\circ$K. The Hamiltonian was given by 
\be
\label{BCHam}
\HH(\s)=
-J
\sum_{\xy} \s_x\s_y
+ 
D \sum_{x\in V} \s_x^2
-
h \sum_{x\in V} \s_x, \qq \s\in\{-1,0,+1\}^V,
\ee
where $J$, $D$, $h$ are constants. The probability of a 
configuration $\s$ was taken proportional to $e^{-\b\HH(\s)}$, $\b=1/(kT)$. 
Capel \cite{CAPEL1,CAPEL2,CAPEL3} 
used molecular field approximations to study the 
ferromagnetic case $J>0$. 
A special case is the system with zero external-field, that is, $h=0$.
For a regular graph with vertex degree $\d$, 
Capel calculated that 
there is a first-order phase transition when
$\frac13 J\d \log 4<D<\frac12 J \d$, and a
second-order phase transition when $D<\frac13 J\d \log 4$. 
For $D> \frac12 J \d$ he predicted that zero states would be dominant.
These non-rigorous results have led to a certain amount
of interest in the\ \BC\ model. According to the physics
literature, there is a first-order transition
even in the low-dimensional setting of $\ZZ^2$. Indeed, in
the phase diagram with parameters $(J,D)$, there is
believed to be a so-called \lq tri-critical point', at which
a line of phase transitions turns from first- to second-order.

The so-called \lq\rc\ representation' of Fortuin and Kasteleyn provides one
of the basic methods for studying Ising and Potts models,
see \cite{MR2014387}--\cite{G-RC} 
and the references therein. Our
target in the current paper is to demonstrate a \rc\ representation
for the \BC\ model with $h=0$. One of the principal advantages
of this approach is that it allows the use of stochastic monotonicity
for the corresponding \rc\ model. Thus, we shall explore
monotonicity and domination methods for the ensuing measure,
and shall deduce some of the structure of the \BC\ model on
$\ZZ^d$.

There is some related literature. A different
approach to a random-cluster representation of the Blume--Capel
model is discussed in \cite{Bouabci-C}, where the target
was to implement a Monte-Carlo method 
of Swendsen--Wang type, \cite{SwendsenWang}.
A related but different problem is the Potts lattice gas, in
which the usual Potts state-space is augmented by an extra site $0$ representing
an empty vertex, see \cite{Al01}. We note the early paper of Hu,
 \cite{HU}, who considered a \rc\ representation 
for the Ising model with general ferromagnetic cell interaction on a 
square lattice. 

For further results on the \BC\ model, see
\cite{BBCK,BCC05,BS89,CO96,FP04,HK02,OM01}.
The usual \rc\ model is summarized in \cite{G02,G-RC}.

The \BC\ model has three local states. There is an extension to
a model with local state space $\{0,1,2,\dots,q\}$ where $q\ge 1$.
We introduce this new model in Section \ref{section:3}, 
where we dub it the \BCP\ (BCP) model. 
We show there how to construct a \rc\ representation of
the BCP model, and we call the corresponding model
the `\drc' (DRC) model. In the BCP model, vertices with state zero do not interact 
further with their neighbours, and the states of the other vertices have 
a Potts distribution. In the \drc\ model, the zero-state
vertices of the BCP model are removed, and the remaining
graph is subject to a conventional \rc\ model.
Note that the \drc\ model is an `annealed' model in the sense
that the dilution is done at random.

The \drc\ model is formulated on a finite graph in Section \ref{section:3}, 
and with boundary conditions on a 
(hyper)cubic lattice in Section \ref{latticeDRC}.
In Section \ref{BoundaryConditionsStochasticOrderings},
we establish stochastic orderings
of measures, and we use these to study phase transitions.
There are two types of stochastic ordering. 
In Section \ref{BoundaryConditionsStochasticOrderings}, we study the
process of vertex-dilution, and we show that the set of remaining vertices
has a law which is both monotonic and satisfies
stochastic orderings with respect to different parameter values.
In Section \ref{EdgeComparisonSection}, we consider the set of open
edges after dilution, and we prove stochastic orderings 
for the law of this set. The results so far are for finite graphs only.

The thermodynamic limit is taken in two steps, 
in Section \ref{section:infiniteVolumeLimit}. 
We prove first the existence of the infinite-volume limit
of the vertex-measure, and the
infinite-volume limits of the full measure and of the BCP measure
follow for $1\le q\le 2$. 
As in the case of the \rc\ model, a certain
amount of uniqueness may be obtained using an argument
of convexity of pressure.
The comparison results for finite graphs carry through
to infinite graphs, and enable a rigorous but incomplete study
of part of the phase diagram of the \BC\ model.
This is summarized in Sections \ref{ptransitions} and \ref{phasediagram},
where it is shown that the rigorous theory of the $q=1$
case gives support for the conjectured phase diagram of the \BC\ model.

\section{Notation}
A finite graph $G=(V,E)$ comprises a vertex-set $V$ and
a set $E$ of edges $e=\la x,y\ra$ having endvertices $x$ and $y$.
We write $x \sim  y$ if $\la x,y\ra \in E$, and we call $x$
and $y$ {\em neighbours} in this case. For simplicity, we
shall assume generally that $G$ has neither loops nor multiple edges.
 The degree $\deg_x$ of a vertex
$x$ is the number of edges incident to $x$.

Let $d\ge 2$.
Let $\ZZ=\{\dots,-1,0,1,\dots\}$, and let
$\ZZ^d$ be the set of all $d$-vectors of integers.
For $x\in \ZZ^d$, we write $x=(x_1,x_2,\dots,x_d)$, and we define
\[
|x| = \sum_{i=1}^d |x_i|.
\]
We write $x \sim y$ if $|x-y|=1$, and we let $\Ed$ be
the set of all unordered pairs $\la x,y\ra$ with $x\sim y$.
The resulting graph $\Ld=(\ZZ^d,\Ed)$ is called the $d$-dimensional
{\em hypercubic lattice}.

Substantial use will be made later of the Kronecker delta,
\[
\d_{u,v} = \begin{cases} 1 &\text{if } u=v,\\
0&\text{if } u\ne v.
\end{cases}
\]

\section{The BCP and DRC measures}\label{section:3}
It is shown in this section how the \BC\ measure on a graph may be 
coupled with a certain `diluted-\rc' measure. 
The \BC\ model has two non-zero local states, labelled
$\pm 1$. Just as in the Ising/Potts case, the corresponding
\rc\ representation is valid for a general number, 
$q$ say, of local states. Therefore, we first define
a `Potts extension' of the \BC\ model with zero external-field.

Let $G=(V,E)$ be a finite graph with
 neither loops nor multiple edges.
Let $q\in\{1,2,3,\dots\}$, and let $\Sigma_q
=\{0,1,2,\dots,q\}^V$. For $\s=(\s_x: x\in V)\in\Si_q$,
we let $E_\s$ be the subset of $E$ comprising all
edges $e=\xy$ with $\s_x\ne 0$, $\s_y\ne 0$.
After a change of notation, the
\BC\ measure with zero external-field 
amounts to the probability measure on $\Si_2$
given by
\be
\pi_2 (\s ) = \frac1 {\ZBC} \exp \left[
-K|E_\s| + 2K\sum_{e \in E} \d_e(\s)
+
\De \sum_{x \in V} \d_{\s_x,0} 
\right],
\qq\s\in\Si_2,\label{BC2def}
\ee
where
\[
\d_e(\s) = \d_{\s_x,\s_y}(1-\d_{\s_x,0}),\qq e=\xy \in E.
\]
Note that
\[
2\d_e(\s) = \s_x\s_y + 1\qq\text{for $e=\la x,y\ra$ and 
$\s_x,\s_y\in\{-1,+1\}$},
\]
and this accounts for the exponent in
(\ref{BC2def}).
The constants $K$ and $\De$ are to be regarded as parameters of the model.
We now define the `\BCP\ (BCP)' probability measure
$\pi_q$ on $\Si_q$ by
\be
\pi_q(\s) = \frac 1 \ZBCP \exp \left[-K|E_\s|+ 2K
\sum_{e \in E} \d_e(\s) + \De \sum_{x \in V} \d_{\s_x,0} 
\right],
\qq\s\in\Si_q, \label{BCPdef}
\ee
where $\ZBCP=\ZBCP_{K,\De,q}$ is the normalizing constant.
We point out that the value $q=1$ is permitted in the above definition
of $\pi_q$.

We turn now to the \rc\ representation of the BCP measure. 
The support of the corresponding \rc-type measure
 is a subset of the product $\Psi\times \Om$ where 
$\Psi = \{ 0,1\} ^V$, and $\Omega=\{ 0,1\}^E$. For 
$\psi=(\psi_x: x\in V) \in\Psi$,
we let 
\[
V_\psi=\{x\in V:\psi_x=1\},\q
E_\psi = \{\xy\in E :x,y\in V_\psi\}. 
\] 
Let $\om=(\om_e: e\in E)\in\Om$. We say that
$\om$ and $\psi$ are {\em compatible\/} if $\om_e=0$
whenever $e\notin E_\psi$, and we write $\Th$
for the set of all compatible pairs $(\psi,\om)\in\Psi\times\Om$. 
Let $\t=(\psi,\om)
\in\Th$. A vertex $x\in V$ is called {\em open} (or $\psi$-open) if
$\psi_x=1$, and is called {\em closed} otherwise.
An edge
$e$ is called {\em open} (or $\om$-open) if $\om_e=1$, and {\em
closed} otherwise. We write $\eta(\om)$ for the set of
$\om$-open edges, and note that $(\psi,\om)\in\Theta$ if
and only if $\eta(\om)\subseteq E_\psi$.
For $\t=(\psi,\om)\in\Th$ and $e\notin E_\psi$,
we say that $e$ has been {\em deleted}. 

Let $\t=(\psi,\om)\in\Th$.
The connected 
components of 
the graph  $(V_\psi,\eta(\om))$ are called {\em open
clusters}, and their cardinality is denoted by $k(\t)$.

The parameters of the \rc\ measure in question are
$a\in(0,1]$, $p\in[0,1)$, $q\in(0,\oo)$, and in addition we write $r=\sqrt{1-p}$.
The {\em \drc\ measure} with parameters $a$, $p$, $q$
is defined to be the probability
measure on $\Psi\times\Om$ given by
\be
\phi (\t) =
\frac 1{\ZDRC}{r^{|E_\psi|}q^{k(\t)}} 
\prod_{x\in V} \left(\frac{a}{1-a}\right)^{\psi_x}
\prod_{e\in E_\psi} \left(\frac{p}{1-p}\right)^{\om_e}
\label{drcdef}
\ee
for  $\t=(\psi,\om)\in\Th$, and $\phi(\t)=0$
otherwise,
where $\ZDRC = \ZDRC_{a,p,q}$ is the normalizing constant. 
The above formula may be interpreted when
$a=1$ as requiring that all vertices be open. 
We note for future use that the projection of $\phi$ onto the 
first component $\Psi$ of the configuration space is the probability measure
satisfying
\be
\Phi(\psi) = \sum_{\om\in\Om} \phi(\psi,\om)
\propto r^{|E_\psi|}\left(\frac a{1-a}\right)^{|V_\psi|} 
\ZRC_{p,q}(V_\psi,E_\psi),
\qq\psi\in\Psi,
\label{projrc}
\ee
where 
\be
\ZRC_{p,q}(W,F)=\sum_{\om\in\{0,1\}^F} q^{k(\om)} 
\left(\frac p{1-p}\right)^{|\eta(\om)|}
\label{rcpdef}
\ee
denotes the partition function of the 
\rc\ model on $G=(W,F)$ with parameters $p$, $q$.
When $F=\es$, we interpret $\ZRC_{p,q}(W,F)$ as $q^{k(W,F)}$,
where $k(W,F)$ is the number of components of the graph.
We speak of $\Phi$ as the `vertex-measure' of $\phi$.

The \drc\ and BCP measures are related to one another
 in very much the same
way as are the \rc\ and Potts measures, see \cite{G-RC}.
This is not quite so obvious as it may first seem,
owing to the factor $r^{|E_\psi|}$ in the definition of $\phi$.
We will not labour the required calculations since they
follow standard routes, but we present the coupling
theorem, and we will summarize some of the necessary facts 
concerning
the conditional measures.

We turn therefore to a coupling between the \drc\ and BCP measures.
Let $\De\in \RR$, $K\in[0,\oo)$, $q\in\{1,2,3,\dots\}$,
and let $a$ and $p$ satisfy
\be
p=1-e^{-2 K},\q \frac{a}{1-a}=
e^{-\De}.
\label{apdef}
\ee
We will define a probability measure $\mu$ on the product
space $\Si_q \times \Psi\times \Om$. This measure $\mu$
will have as support the subset
$\sS\subseteq \Sigma_q\times\Psi
\times\Om$ comprising all triples $(\s,\psi,\om)$ such that:
\begin{romlist}
\item $(\psi,\om)\in\Th$,
\item $\psi_x= 1-\d_{\s_x,0}$ for all $x\in V$, that is, $\psi_x=0$ if and only if $\s_x= 0$, and
\item for all $e=\xy\in E$, if $\s_x\ne\s_y$ then $\om_e=0$.
\end{romlist}
We define $\mu$ by
\[
\mu(\s, \psi,\om ) = \begin{cases}
\displaystyle\frac 1{Z}{r^{|E_\psi|}} 
\prod_{x\in V} \left(\frac{a}{1-a}\right)^{\psi_x}
\prod_{e\in E_\psi} \left(\frac{p}{1-p}\right)^{\om_e}
&\text{if } (\s,\psi,\om)\in\sS,\\
0&\text{otherwise}.
\end{cases}
\]

\begin{thm}
\label{thm:couple-BC-DRC}
Let $q\in\{1,2,3,\dots\}$, let $\De\in \RR$, $K\in[0,\oo)$
 and let $a$, $p$ satisfy
\eqref{apdef}.
The marginal measures of $\mu$ on $\Sigma_q$ 
and on $\Psi\times\Omega$, respectively,
are the BCP and \drc\ measures with respective parameters 
$K$, $\De$, $q$ and
$a$, $p$, $q$.
\end{thm}

\begin{proof}  
Let $\s \in \Sigma_q$. We fix $\psi$ by  $\psi_x=1-\d_{\s_x,0}$ for
all $x\in V$, so that $E_\psi = E_\s$. By (\ref{apdef}),
\[
\prod_{x\in V} \left(\frac{a}{1-a}\right)^{\psi_x}=
\exp\left[- \De \sum_{x \in V} (1-\d_{\s_x,0})\right].
\]
By summing over all $\om$ such that $(\s,\psi,\om)\in\sS$,
\begin{align*}
\sum_{\om}  r^{|E_\psi|} \prod_{e\in E_\psi}
\Big(\frac{p}{1-p}\Big)^{\om_e}
&=
r^{|E_\psi|} \prod_{e\in E_\psi} \left[ 1+
\Big(\frac{p}{1-p}\Big)\d_e(\s) \right]\\
&=\exp\left[-K|E_\psi| + 2K\sum_{e\in E} \d_e(\s) \right].\\
\end{align*}
By (\ref{BCPdef}),
\[
\sum_{(\psi,\om)\in\Th} \mu(\s,\psi,\om) \propto \pi_q(\s),\qq \s\in\Si_q.
\]
Equality must hold here, since each side is a probability
mass function.  This proves that the marginal of $\mu$ on $\Si_q$
is indeed the BCP measure $\pi_q$. 

Turning to the second marginal, we fix 
$\t=(\psi,\om)\in\Th$,
and let $\sS(\t)$ be the set of all $\s\in\Si_q$
such that $(\s,\psi,\om)\in\sS$. We have that $\s_x=0$
if and only if $\psi_x=0$. The only further constraint
on $\s$ is that it is constant on each cluster of
$(V_\psi,\eta(\om))$.
There are $k(\t)$ such clusters,
and therefore $|\sS(\t)|=q^{k(\t)}$.
It follows that
\[
\sum_{\s} \mu(\s,\psi,\om)= \phi(\t),\qq 
\t=(\psi,\om)\in \Th,
\]
as required.
\end{proof}

We make some observations based on Theorem \ref{thm:couple-BC-DRC} 
and the method of proof.
First, subject to (\ref{apdef}),
\be
\ZBCP_{K,\De,q}=\ZDRC_{a,p,q} \cdot e^{|V|\De}.  \label{eqpartfns}
\ee
Secondly, the conditional measure of $\mu$, given the
pair $(\psi,\om)\in \Psi\times\Om$, is that obtained as follows:
\begin{letlist}
\item for $x\in V$, $\s_x=0$ if and only if $\psi_x=0$,
\item the spins are constant on every cluster of the graph
$(V_\psi,\eta(\om))$, and each such spin is uniformly distributed
on the set $\{1,2,\dots,q\}$,
\item the spins on different clusters are independent
random  variables.
\end{letlist}
Thirdly, the conditional measure of $\mu$, given the
spin vector $\s\in\Sigma_q$, is that obtained as follows:
\begin{romlist}
\item for $x\in V$, $\psi_x = 0$ if and only if $\s_x=0$,
\item $(\psi,\om)\in\Th$,
\item the random variables $(\om_e: e\in E_\psi)$ are independent,
\item for $e=\la x,y\ra \in E_\psi$, $\om_e=0$ if $\s_x\ne \s_y$, and
$\om_e=1$ with probability $p$ if $\s_x=\s_y$.
\end{romlist} 
In particular, conditional on the set
$\{x\in V: \s_x=0\}$, the joint distribution of $\s$ and $\om$ is the
usual coupling of the Potts and \rc\ measures 
on the graph $G_\psi=(V_\psi,E_\psi)$.

As two-point correlation function in the BCP model, we may take
the function 
\be\label{twoptcf}
\tau_q(x,y)=\pi_q(\s_x=\s_y\ne 0)- \frac 1q \pi_q(\s_x\s_y\ne 0),
\qq x,y\in V. 
\ee
This is related as follows to the two-point connectivity function 
of the \drc\ model.
For $x,y\in V$, we write $x\lra y$ if there
exists a path of $\om$-open edges joining $x$ to $y$.
Similarly, for $A,B\subseteq V$, we write $A \lra B$ if there
exist $a\in A$ and $b\in B$ such that $a\lra b$.

\begin{thm}\label{corrconn}
Let $\De\in\RR$, $K\in[0,\oo)$, $q\in\{1,2,\dots\}$, and let
$a$, $p$ satisfy \eqref{apdef}.
The corresponding \drc\ measure $\phi$ and BCP measure $\pi_q$ on
the finite graph $G=(V,E)$ are such that
$$
\tau_q(x,y) = 
(1-q^{-1}) \phi(x\lra y),
\qq x,y\in V.
$$
\end{thm}

The proof follows exactly that of the corresponding statement
for the \rc\ model, see for example \cite{G-RC}.

Two particular values of $q$ are special, namely $q=1,2$. From the above,
the \drc\ measure with $q=2$ corresponds to
the \BC\ measure.
Theorem \ref{thm:couple-BC-DRC} is valid
with $q=1$ also. The BCP model with $q=1$ has two local states
labelled $0$ and $1$. By (\ref{BCPdef}), 
the Hamiltonian may be written as
\begin{align*}
\HH_q(\s) &= K|E_\s|  - 2K\sum_{e\in E} \d_e(\s) -\De\sum_{x\in V}\d_{\s_x,0}\\
&= -\De |V| -K|E_\s| + \De\sum_{x\in V} \s_x\\
&= -\De |V| -K\sum_{e=\xy\in E} \s_x\s_y + \De \sum_{x\in V} \s_x,
\qq\s\in\Si_1.
\end{align*}
We make the change of variables $\eta_x=2\s_x-1$, to find that
\[
\HH_q(\s) =-\tfrac12 \De|V|  -\tfrac 14 K|E| - 
J \sum_{e=\xy\in E} \eta_x\eta_y 
-\sum_{x\in V} h_x \eta_x,
\]
where $J=\frac14 K$ and $h_x=\frac14(K\deg_x - 2\De)$.
That is, we may work with the altered Hamiltonian
\be\label{Isingh}
\HH'_q(\s) =
-J \sum_{e=\xy\in E} \eta_x\eta_y         
-\sum_{x\in V} h_x \eta_x,
\ee
which is recognised as that of the Ising model with
edge-interaction $J$ and `local' external field $(h_x:x\in V)$.
If $G$ is regular with (constant)
vertex-degree $\d$, then $h_x=h=\frac14(K\d -2\De)$ for all $x\in V$. 
That is,
the BCP model with $q=1$ is, after a re-labelling
of the local states $0$, $1$, an Ising model
with edge-interaction $J$ and external field 
$h$.  A great deal is known about this model, 
and we shall make use
of this observation later.

\section{The lattice DRC model}\label{latticeDRC}
Until further notice, we shall study the \drc\ model
rather than the BCP model, and thus we 
take $q$ to be a positive real (number).
The model has so far been defined on a finite graph only.
In order to pass in Section \ref{section:infiniteVolumeLimit}
to the infinite-volume limit on 
$\Ld$, we shall next introduce the concept of boundary conditions.

Let $V$ be a finite subset of $\ZZ^d$, and let 
$E$ be the subset of $\Ed$ comprising all edges having at least one
endvertex in $V$.  We write $\La=(V,E)$, 
noting that $\La$
is not a graph since it contains edges adjacent 
to vertices outside $V$. Any such $\La$ is
called a {\it region\/}. The corresponding graph $\La^+
=(V^+,E)$ is defined as the subgraph of $\Ld$ induced by $E$.
We write $\pd\La=V^+\sm V$.
The lattice $\Ld$ is regular
with degree $\d=2d$.

Let $\Psi=\{0,1\}^{\ZZ^d}$ and $\Om=\{0,1\}^{\Ed}$.
Let $\Th$ be the set of compatible vertex/edge configurations 
$(\psi,\om)\in
\Psi\times\Om$ satisfying $\eta(\om)\subseteq \Ed_\psi$.
Each $\l=(\kappa,\rho)\in\Th$ may be viewed as a boundary condition on the region $\La$,
as follows.
Let $\Th_\La^\l$ be the subset of $\Th$ containing 
configurations that agree with $\l$ on $\Ld\sm \La$,
in that $\Th_\La^\l$ contains all $(\psi,\om)$ with
$\psi_x=\kappa_x$ for $x\notin V$, $\om_e=\rho_e$ for $e\notin E$.
Let $\phi^\l_{\La,a,p,q}$ denote the \drc\ measure on
$\La$ with boundary condition $\l$, that is,
\be
\label{new34}
\phi^\l_{\La,a,p,q} (\t) =
\frac 1{\ZDRC}{r^{|E_\psi|}q^{k(\t,\La)}} 
\prod_{x\in V} \left(\frac{a}{1-a}\right)^{\psi_x}
\prod_{e\in E_\psi} \left(\frac{p}{1-p}\right)^{\om_e},
\ee
if $\t=(\po)\in\Th^\l_{\La}$, and
$\phi^\l_{\La,a,p,q} (\t) =0$ otherwise.
Here, $E_\psi=\{\langle x,y\rangle\in E:
\psi_x=\psi_y=1\}$, 
$k(\t,\La)$ is the number of
open clusters of $(\ZZ^d_\psi,\eta(\om))$  that 
intersect $V^+$ , and
$\ZDRC = \ZDRC_{\La,\l,a,p,q}$ is a normalizing constant.
See \eqref{drcdef}, and recall that $r=\sqrt{1-p}$.

The probability measure $\phi^\l_{\La,a,p,q}$ is
supported effectively on the product $\Psi_V\times
\Om_E$ where $\Psi_V=\{0,1\}^V$ and
$\Om_E=\{0,1\}^E$.  We write $\Phi^\l_{\La,a,p,q}$ for 
its marginal (or `projected') measure 
on the first coordinate
$\{0,1\}^V$ of this space, given as follows.
Let $\l=(\kappa,\rho)\in\Th$, let $\La=(V,E)$ be a region,
and let $\Psi_\La^\l$ be the set of all $\psi\in\Psi$ that agree
with $\kappa$ off $V$.
For $\psi\in\Psi$, let $\La(\psi)$ denote the subgraph of $\La^+$
induced by the $\psi$-open vertices. Suppose 
$\psi\in\Psi_\La^\l$.
Let 
$\ZRC_{\l,p,q}(\La(\psi))$ denote the partition function
of the \rc\ model on $\La(\psi)$ with boundary condition $\l$,
see \eqref{rcpdef}.
(This boundary condition is to be interpreted as: 
two vertices $u,v\in V^+$ are deemed to be connected off $\La^+$
if there exists a path from $u$ to $v$ of $\rho$-open
edges of $\Ed\sm E$.) As in \eqref{projrc},
\begin{align}
\Phi^\l_{\La,a,p,q}(\psi) &= \sum_{\om\in\Om_E} 
\phi^\l_{\La,a,p,q}(\psi,\om)\nonumber\\
&\propto r^{|E_\psi|}\left(\frac a{1-a}\right)^{|V_\psi|} 
\ZRC_{\l,p,q}(\La(\psi)),
\label{projrc2}
\end{align}
for $\psi\in\Psi_V$,
where $V_\psi = \{v\in V: \psi_v=1\}$. There is a slight abuse
of notation here, in that $\psi$ has been used as a member of 
both $\Psi$ and $\Psi_V$.

Two especially interesting situations arise when $p=0$ and/or $q=1$.

\noindent
(a) {\em Product measure.} 
If $p=0$ then $\phi^\l_{\La,a,p,q}$ is a product
measure, and may therefore be extended to a product measure
$\phi_{a,0,q}$ on $\Ld$ under which each
vertex is open
with probability $qa/(1-a + qa)$, and
each edge is almost-surely closed.
There exists,
$\phi_{a,0,q}$-almost-surely, an infinite open
vertex-cluster (\resp, infinite closed vertex-cluster) if
$qa/(1-a + qa)>\pcsite$ (\resp, $(1-a)/(1-a + qa)>\pcsite$),
where $\pcsite$ denotes the critical probability of site 
percolation on $\Ld$.

\noindent
(b) {\em Ising model with external field.} 
Let 
$q=1$, and recall from the end of Section \ref{section:3}
that the BCP model is essentially an Ising model
with edge-interaction $J=\frac14 K$ and local
external field $h_x=\frac14(K\deg_x-2\De)$.
For the sake of illustration, consider the
box $B_n=[-n,n]^d$ of $\LL^d$ with periodic boundary conditions,
so that $\deg_x=\d=2d$ for all $x$. Then
\be\label{new36}
J=
-\tfrac18 \log(1-p),\q
h=\tfrac12(Kd-\De)=\tfrac12\log\left(\frac a{(1-a)(1-p)^{d/2}}\right).
\ee

On passing to the limit as $n\to\oo$, we obtain
an infinite-volume Ising model with parameters $J$, $h$.
If we restrict ourselves to pairs $a$, $p$ such that
$h=0$, there is a critical value 
$\Kc(d)$ of $K$ given by
$\Kc(d)=-2\log(1-\pic)$ where $\pic=\pic(d)$ is
the critical edge-parameter of the \rc\ model on $\LL^d$
with cluster-weighting parameter $2$.
Rewritten in terms of $a$ and $p$,
the phase diagram possesses a special point
$(\ol a,\ol p)$, where
\be\label{olaolp}
\ol a= \frac{(1-\pic)^{2d}}{1+(1-\pic)^{2d}},
\q \ol p=1-(1-\pic)^4.
\ee
By a consideration of the associated \rc\ measure or otherwise, we
deduce that there is a line of first-order phase transitions along the arc
\be\label{new37}
\frac a{1-a} = (1-p)^{d/2},\qq \ol p <p<1.
\ee
To the left (\resp, right)
of this arc in $(a,p)$ space (see Figure \ref{newfig3}
for the case $d=2$), there is an infinite cluster of $0$-state
(\resp, $1$-state)
vertices. As the arc is crossed from left to right, there is
a discontinuous increase in the density of the infinite
$1$-state cluster. Related issues concerning the percolation
of $\pm$-state clusters in the zero-field Ising model 
are considered in \cite{ABL}.

We note when $d=2$ that $\pic(2)=\sqrt 2/(1+\sqrt 2)$, so that
\be\label{olaolp2}
\ol a = \frac1{1+(1+\sqrt 2)^4},\q
\ol p = 1-(1+\sqrt 2)^{-4}.
\ee

\section{Stochastic orderings of vertex-measures}
\label{BoundaryConditionsStochasticOrderings}
Many of the results of this section have 
equivalents for general finite
graphs, but we concentrate here on subgraphs
of the lattice $\Ld=(\ZZ^d,\Ed)$.
While the route followed here is fairly standard,
some of the calculations are novel.
The vertex-measure
$\Phi^\l_{\La,a,p,q}$ plays an important part in the stochastic
orderings relevant to the BCP model, and we turn next to
its properties, beginning with a reminder about
orderings.

Let $I$ be a finite set, and let $\Sigma=\{0,1\}^I$ be viewed as a 
partially ordered set. 
For $J\subseteq I$ and $\s\in\Sigma$,
we write $\s^J$ for the configuration that 
equals 1 on $J$ and
agrees with $\s$ off $J$. 
If $J=\{i\}$ or $J=\{i,j\}$ we may
abuse notation by removing the braces.
Let $\mu_1$, $\mu_2$
be probability measures on $\Sigma$. We write $\mu_1\lest\mu_2$,
and say that $\mu_1$ is {\it stochastically dominated\/} by 
$\mu_2$,
if $\mu_1(f) \le \mu_2(f)$ for all increasing
functions $f:\Sigma \to \RR$. A probability measure $\mu$
on $\Sigma$ is said to be strictly positive if
$\mu(\s)>0$ for all $\s\in\Sigma$. If $\mu_1$, $\mu_2$
are strictly positive, then $\mu_1\lest \mu_2$ if the pair satisfies
the so-called Holley condition,
\begin{equation}
\mu_2(\s_1\vee\s_2)\mu_1(\s_1\wedge\s_2) 
\ge \mu_1(\s_1)\mu_2(\s_2),
\qq \s_1,\s_2\in\Sigma.
\label{Holleycond}
\end{equation}
Here, $\vee$ denotes the coordinatewise maximum, 
and $\wedge$ the
coordinatewise minimum. It is standard (see \cite{G-RC}, Section 2.1)
that it suffices to check \eqref{Holleycond} for pairs
of the form $(\s_1,\s_2) = (\s^i,\s)$ and 
$(\s_1,\s_2)=(\s^i,\s^j)$, for $\s\in\Sigma$ and $i,j\in I$.

A probability measure
$\mu$ on $\Sigma$ is said to be
{\it positively associated\/} if
\[
\mu(A \cap B) \ge \mu(A)\mu(B)
\]
for all increasing events $A,B\subseteq \Sigma$.
For $\tau\in\Sigma$ and $J\subseteq I$, 
let $\Sigma_J^\tau$ be the
subset of $\Sigma$ containing all $\s\in\Sigma$ with
$\s_i=\tau_i$ for $i\notin J$.
The measure $\mu$ is said to be 
{\it strongly positively associated\/} if, for all pairs $\tau$, $J$, the conditional
measure, given $\Sigma_J^\tau$, is positively associated
when viewed as a measure on $\{0,1\}^J$. 
The measure $\mu$
is called {\it monotonic\/} if, for all $i\in I$,
$\mu(\s_i=1\mid \Sigma_i^\tau)$ is a non-decreasing 
function of $\tau$. It is standard
(see \cite{G-RC}, Section 2.2) that a strictly positive
probability measure $\mu$ on $\Sigma$ is strongly positively
associated (\resp, monotonic) if and only if
it satisfies the
so-called FKG condition:
\begin{equation}
\mu(\s_1\vee\s_2)\mu(\s_1\wedge\s_2) 
\ge \mu(\s_1)\mu(\s_2),
\qq \s_1,\s_2\in\Sigma.
\label{FKGcond}
\end{equation}
Furthermore,
it suffices to check \eqref{FKGcond} for pairs
of the form
$(\s_1,\s_2)=(\s^i,\s^j)$, for $\s\in\Sigma$ and $i,j\in I$.
Further discussions of the FKG and Holley inequalities
may be found in
\cite{MR0309498,G-RC,MR0341552}. 

The proofs of the following theorems will be found later in this section.

\begin{thm}
\label{thm:positively-associated}
Let $\La=(V,E)$ be a region, let $\l\in\Th$, and let 
$a\in(0,1),p\in[0,1)$. The probability measure
$\Phi^\l_{\La,a,p,q}$ is strongly positively associated, and hence monotonic, if $q\in[1,2]$.
\end{thm}

The condition $q\in[1,2]$ is important.
If $q>2$, then
strong positive-association does not hold for all $p\in(0,1)$.
The conclusion would be similarly false
for the full \drc\ measure 
even for $q\in[1,2]$.
For example, let $G$ be the graph with exactly two vertices $x$, $y$
joined by a single edge $e$, and consider the associated measure
$\phi_{a,p,q}$ with $a,p\in(0,1)$ and $q\in(0,\oo)$. Then, with $r=\sqrt{1-p}$,
\begin{align*}
\phi(\psi_y=1 \mid \psi_x=0,\ \om_e=0)=&\frac{qa}{qa+1-a},\\
\phi(\psi_y=1 \mid \psi_x=1,\ \om_e=0)=&\frac{qar}{qar+1-a}.
\end{align*}
The first term exceeds the second strictly, and hence 
$\phi_{a,p,q}$ is not monotone on the product space $\{0,1\}^V
\times\{0,1\}^E$.

We prove next that $\Phi^\l_{\La,a,p,q}$ is
increasing in $\l$, so long as $q\in[1,2]$.

\begin{thm}
\label{thm:ordering}
Let $\La=(V,E)$ be a region, and let
$a\in(0,1),p\in[0,1)$ and $q\in [1,2]$.
If $\laa\leq\lab$ then $\Phi^\laa_{\La,a,p,q}
\lest \Phi^\lab_{\La,a,p,q}$.
\end{thm}

The two theorems above will be proved by checking certain
inequalities related to (\ref{Holleycond}) and (\ref{FKGcond}).
It is convenient
to make use of a subsidiary proposition for this,
and we state this next, beginning with
some notation. For a region $\La=(V,E)$, we
abbreviate to $\Phi_i$ the marginal 
(or projected) measure 
on the space $\Psi_V$ of the \drc\ measure
$\phi^{\lai}_{\La,a_i,p_i,q_i}$.
We abbreviate to $\mu_{\La,\psi}^i$ the usual
\rc\ measure on $\La(\psi)$ with boundary condition
$\lai$ and parameters $(p_i,q_i)$. 
For $w\in \ZZ^d$, let $I_w\subseteq \Om$
be the event that $w$ has no incident $\om$-open edges.

\begin{prop}\label{thm:big-ugly-lemma}
Let $\laa,\lab\in\Th$,
$a_i\in(0,1)$, $p_i\in[0,1)$ for $i=1,2$, and
$q_1\in[1,\oo)$, $q_2\in[1,2]$.
Let $\psi\in\Psi$, let $\La=(V,E)$ be a region, and let $x\in V$ 
be such that $\psi_x=0$. 
Let $b=b(x,\psi)$ denote
the number of edges of $E$ of the form
$\langle x,z\rangle$
with $\psi_z=1$. If 
\begin{equation}
\label{eqn:tech-lemma-start}
q_2\left(\frac{a_2}{1-a_2}\right)
\frac{(1-p_2)^{b/2}}{\mu_{\La,\psi^x}^2(I_x)}
\ge 
q_1\left(\frac{a_1}{1-a_1}\right)
\frac{(1-p_1)^{b/2}}{\mu_{\La,\psi^x}^1 (I_x)},
\end{equation}
then
\begin{align}\label{eqn:tech-lemma-statement1}
\Phi_2(\psi^{x})\Phi_1(\psi) &\ge
\Phi_1(\psi^x)\Phi_2(\psi),\\
\Phi_2(\psi^{x,y})\Phi_1(\psi) &\ge
\Phi_1(\psi^x)\Phi_2(\psi^y),\qq y\in V\sm V_\psi,\ y\ne x.
\label{eqn:tech-lemma-statement2}
\end{align} 
\end{prop}

We examine next the monotonicity properties of 
$\Phi^\l_{\La,a,p,q}$ as $a$, $p$, $q$ vary.
Recall that $\d=2d$.

\begin{thm}
\label{thm:comparison-inequality}
Let $\La=(V,E)$ be a region, and let $\l\in\Th$.
Let $a_i\in (0,1)$, $p_i\in[0,1)$, and $q_i\in[1,2]$ for $i=1,2$,
and let $\Phi_i$ be as above.
Each of the following is a sufficient condition for the 
stochastic inequality $\Phi_1 \lest \Phi_2${\rm:}
\begin{romlist}
\item
that $a_1\le a_2$, $p_1\le p_2$, and $q_1=q_2$,
\item
that 
\[
q_2\left(\frac{a_2}{1-a_2}\right)
\ge
q_1\left(\frac{a_1}{1-a_1}\right)
(1-p_1)^{-\d/2},
\]
\item
that $p_1\le p_2$, 
$q_1 \ge q_2$,
and
\begin{equation}
q_2\left(\frac{a_2}{1-a_2}\right) (1-p_2)^{\d/ 2}
\ge
q_1\left(\frac{a_1}{1-a_1}\right) (1-p_1)^{\d/ 2}, 
\label{compcond}
\end{equation}
\item 
that
$q_1\leq q_2$, \eqref{compcond} holds, and
\[
\frac{p_2}{q_2(1-p_2)}\geq \frac{p_1}{q_1(1-p_1)}.
\]
\end{romlist}
\end{thm}

In the next section we shall pass to 
infinite-volume limits
along increasing sequences of regions. In preparation for 
this, we note two further properties of stochastic 
monotonicity. The two extremal boundary conditions are
the vectors $\bzero = (0,0) \in \Psi\times \Om$ and
$\bone = (1,1)\in\Psi\times\Om$.

\begin{thm}
\label{thm:monotone-limits}
Let $a\in(0,1),p\in[0,1)$, $q\in[1,2]$, and
let $\La_1$, $\La_2$ be regions with
$\La_1\subseteq \La_2$. Then 
\[
\Phi^\bzero_{\La_1,a,p,q}\lest
\Phi^\bzero_{\La_2,a,p,q},\qq 
\Phi^\bone_{\La_1,a,p,q}\gest
\Phi^\bone_{\La_2,a,p,q}.
\]
\end{thm}

It is noted that the boundary conditions
$b=\bzero, \bone$ contain information concerning both vertex
and edge configuration off $\La$. By \eqref{projrc2},
only the external edge configuration is in fact relevant.
The above inequalities for the vertex-measures $\Phi_{a,p,q}$ 
imply a degree of monotonicity of the full \drc\ measure $\phi_{a,p,q}$.
We shall not explore this in depth, but restrict ourselves to two facts
for later use.

\begin{thm}\label{phimon}
Let $a\in(0,1)$, $p\in [0,1)$, $q\in[1,2]$, and $\l\in\Theta$. 
For any region $\La$, the \drc\ measure
$\phi_{\La,a,p,q}^\l$ is stochastically non-decreasing in $a$, $p$, and $\l$.
\end{thm}

A probability measure on a product space $\{0,1\}^I$ is
said to have the finite-energy property if, 
for all $i\in I$, 
the law of the state of $i$, conditional on the states of all other indices,
is (almost surely) strictly positive. See \cite{G-RC}.

\begin{thm}
\label{thm:finite-energy-Phi} Let $a\in(0,1)$, $p\in[0,1)$,
$q\in[1,2]$, 
$\l\in\Th$, and
let $\La$ be a region. The probability measure
$\Phi^\l_{\La,a,p,q}$ has the 
finite-energy property, and indeed,
\[
\frac{qa}{1-a+qa} \leq
\Phi^\l_{\La,a,p,q}(J_x \mid \sT_x)
\leq \frac{  a q}{ a q + (1-a)r^\d},
\qq \Phi^\l_{\La,a,p,q}\text{-a.s.}, 
\]
where $J_x\subseteq \Psi$ is the event that $x$ is
open, and $\sT_x$ is the $\s$-field of $\Psi$ generated by the states of
vertices other than $x$.
\end{thm}
We turn now to the proofs, and begin with a 
lemma. 

\begin{lem}
\label{thm:remove-vertex}
Under the conditions of Proposition 
\ref{thm:big-ugly-lemma}, and with $x,y\in V\sm V_\psi$,
\[ 
\mu^2_{\La,\psi^{x,y}}(I_x) \leq
\mu^2_{\La,\psi^x}(I_x) r_2^f,
\]
where $r_2=\sqrt{1-p_2}$ and $f\in\{0,1\}$ is the number of edges 
of\/ $\Ld$ with endvertices $x$, $y$.
\end{lem}

\begin{proof}
We note the elementary inequality
\begin{equation}
\frac{q(1-p)}{p+q(1-p)} \le \sqrt{1-p},\qq p\in[0,1],\ q\in[1,2].
\label{element}
\end{equation}
Let $B$ (\resp, $C$) be the set of $b$ (\resp, $c$) edges 
joining  $x$ (\resp, $y$) to $\psi$-open 
vertices of $V^+$, and let $F$ be the set of edges with 
endvertices $x$, $y$.
Let $B_0$ (\resp, $C_0$, $F_0$) be the (decreasing) event that all 
edges in $B$ (\resp, $C$, $F$) are closed. Since
a \rc\ measure with $q\ge 1$ is positively associated, 
\[
\mu_{\La,\psi^{x,y}}^2(I_x)
\leq \mu_{\La,\psi^{x,y}}^2 (B_0 \cap F_0 \mid C_0).
\]
By an elementary property of \rc\ measures,
see \cite{G-RC},
\begin{align*}
\mu^2_{\La,\psi^{x,y}} (B_0 \cap F_0 \mid C_0)
&= \mu^2_{\La\sm C,\psi^{x,y}} (B_0 \cap F_0)\\
&= \mu^2_{\La\sm C,\psi^{x,y}} (B_0\mid F_0)
\mu^2_{\La\sm C,\psi^{x,y}}(F_0),
\end{align*}
where $\La\sm C$ is obtained from $\La$ by deleting all edges in $C$.
In $\La(\psi^{x,y})\sm C$, the only possible neighbour 
of $y$ is $x$, whence,
for $f=|F|=0,1$,
\[
\mu^2_{\La\sm C,\psi^{x,y}}(F_0)
= \frac{q_2(1-p_2)^f}{1+(q_2-1)(1-p_2)^f}
\leq (1-p_2)^{f/2}=r_2^f,
\]
where we have used \eqref{element} and the fact that
$q_2\leq 2$.
Similarly,
\[
\mu^2_{\La\sm C,\psi^{x,y}} (B_0\mid F_0)=
\mu^2_{\La,\psi^{x}}(B_0)
=\mu^2_{\La,\psi^{x}}(I_x),
\]
and the claim follows.
\end{proof}

\begin{proof}[Proof of Proposition \ref{thm:big-ugly-lemma}]
We prove \eqref{eqn:tech-lemma-statement2} only,
the proof of \eqref{eqn:tech-lemma-statement1} is similar and simpler.
Inequality \eqref{eqn:tech-lemma-start} implies
by Lemma \ref{thm:remove-vertex} that
\begin{equation}
\label{eqn:tech-lemma1}
q_2\left(\frac{a_2}{1-a_2}\right)\frac{r_2^{b+f}}{\mu^2_{\La,\psi^{x,y}}(I_x)}
\ge
q_1\left(\frac{a_1}{1-a_1}\right)\frac{r_1^{b}}{\mu^1_{\La,\psi^{x}} (I_x)} ,
\end{equation}
where $f$ is the number of edges of $\LL^d$ joining $x$ and $y$.
Let $\ZRC_{\l,p,q}(G)$ be the partition function of
the \rc\ model on a graph $G$ with parameters $p$, $q$ and
boundary condition $\l$, see \eqref{rcpdef}.
We have that
\begin{equation*}
\mu^1_{\La,\psi^{x}}(I_x) = q_1
\frac{\ZRC_{\laa,p_1,q_1}(\La(\psi))}
{\ZRC_{\laa,p_1,q_1}(\La(\psi^{x}))},\qq
\mu^2_{\La,\psi^{x,y}}(I_x)=
q_2\frac{\ZRC_{\lab,p_2,q_2}(\La(\psi^{y}))}
{\ZRC_{\lab,p_2,q_2}(\La(\psi^{x,y}))}.
\end{equation*}
We substitute these into \eqref{eqn:tech-lemma1} to find that
\begin{multline*}
\left(\frac{a_2}{1-a_2}\right)
\ZRC_{\lab,p_2,q_2}(\La(\psi^{x,y}))
\ZRC_{\laa,p_1,q_1}(\La(\psi))
r_2^{b+f}\\
\ge
\left(\frac{a_1}{1-a_1}\right)
\ZRC_{\laa,p_1,q_1}(\La(\psi^x))
\ZRC_{\lab,p_2,q_2}(\La(\psi^y))
r_1^{b}.
\end{multline*}
Now, $|V(\psi^x)\sm V(\psi)|=1$ and $|E(\psi^x)\sm E(\psi)|=b$ 
where $V(\psi) = V\cap \ZZ^d_\psi$ and 
$E(\psi)=E\cap \Ed_\psi$, so that
\begin{align*}
&\left(\frac{a_2}{1-a_2}\right)^{|V(\psi^{x,y})|}
\ZRC_{\lab,p_2,q_2}(\La(\psi^{x,y}))
r_2^{|E(\psi^{x,y})|}\\
&\hskip4cm \times \left(\frac{a_1}{1-a_1}\right)^{|V(\psi)|}
\ZRC_{\laa,p_1,q_1}(\La(\psi))r_1^{|E(\psi)|}\\
&\hskip1cm \ge \left(\frac{a_1}{1-a_1}\right)^{|V(\psi^x)|}
 \ZRC_{\laa,p_1,q_1}(\La(\psi^x))r_1^{|E(\psi^x)|}\\
&\hskip4cm\times
\left(\frac{a_2}{1-a_2}\right)^{|V(\psi^y)|}
 \ZRC_{\lab,p_2,q_2}(\La(\psi^y))r_2^{|E(\psi^y)|}.
\end{align*}
As in \eqref{projrc2}, 
\begin{equation*}
\Phi_i(\psi)=\sum_{\om\in\Om_E} \phi^\lai_{\La,a_i,p_i,q_i}(\po)
\propto r_i^{|E(\psi)|}  
\left(\frac{a_i}{1-a_i}\right)^{|V(\psi)|}
\ZRC_{\l_i,p_i,q_i}(\La(\psi)),\\
\end{equation*}
for $\psi \in\Psi_V$, and \eqref{eqn:tech-lemma-statement2} follows.
\end{proof}

\begin{proof}[Proof of Theorem \ref{thm:positively-associated}]
We apply Proposition \ref{thm:big-ugly-lemma} with
$a_i=a$, $p_i=p$, $q_i=q$, and $\l_i=\l$. Inequality
\eqref{eqn:tech-lemma-start} is a triviality since $\mu^1_{\La,\psi}
= \mu^2_{\La,\psi}$ for every $\psi$. By
\eqref{eqn:tech-lemma-statement2} and the comment after
\eqref{FKGcond}, 
$\Phi^\l_{\La,a,p,q}$ satisfies the FKG condition \eqref{FKGcond}, and the claim follows.
\end{proof}

\begin{proof}[Proof of Theorem \ref{thm:ordering}]
Since $\l_1\le \l_2$, $\mu^1_{\La,\psi}
\lest \mu^2_{\La,\psi}$ for every $\psi\in\Psi$. Now, $I_x$
is a decreasing event, whence $\mu^1_{\La,\psi}(I_x)
\ge \mu^2_{\La,\psi}(I_x)$. By Proposition
\ref{thm:big-ugly-lemma} and the comment after \eqref{Holleycond},
the $\Phi_i=\Phi_{\La,a,p,q}^{\l_i}$
satisfy the Holley condition \eqref{Holleycond}, and the claim
follows.
\end{proof}

\begin{proof}[Proof of Theorem 
\ref{thm:comparison-inequality}]
In each case, we shall apply Proposition \ref{thm:big-ugly-lemma} 
and appeal to the Holley condition
\eqref{Holleycond} and the comment thereafter. 
It suffices to check
\eqref{eqn:tech-lemma-start} for every relevant vertex $x$.
We recall some basic facts about \rc\ measures 
to be found in, for example, \cite{G-RC}. Let $G=(W,F)$
be a graph and let $\mupq$ be the \rc\ measure
on $\{0,1\}^F$ with parameters $p\in[0,1]$, $q\in[1,\oo)$. 
By the comparison inequalities,
\begin{equation}
\frac p{p+q(1-p)} \le
\mupq(\text{$f$ is open}) \le p,\qq f\in F,
\label{edgedens}
\end{equation}
and, if $x\in W$ has degree $b$,
\begin{equation}
(1-p)^b \le \mupq(I_x) 
\le \left( 1- \frac {p}{p+q(1-p)} \right)^b.
\label{rccomp}
\end{equation}
We note from \eqref{element} that
\begin{equation}
\left( 1- \frac {p}{p+q(1-p)} \right)^b
\le (1-p)^{b/2},\qq p\in[0,1],\  q\in[1,2].
\label{element2}
\end{equation}

\noindent(i):
We may adapt the exponential-steepness argument of
\cite{MR1480028}, as in Section 2.5 of \cite{G-RC}, to the decreasing 
event $I_x$ to obtain, in the above
notation,
\begin{equation}
\frac d{dp}\log\mupq(I_x) \le -\frac1 {p(1-p)}
\sum_{f:\, f\sim x} \mupq(\text{$f$ is open}),
\label{elem3}
\end{equation}
where the sum is over the $b$ edges $f$ with endvertex $x$.
Let $q\in [1,2]$.
By \eqref{edgedens},
\[
\frac d{dp}\log\mupq(I_x) \le -\frac1{p(1-p)} \sum_{f:\, f\sim x}\frac p{p+q(1-p)}
\le - \frac b{2(1-p)}.
\]
We integrate from $p_1$ to $p_2$ and apply to the measures $\mu_i^{\La(\psi)}$
to obtain that
\[
\frac {\mu^2_{\La,\psi}(I_x)} {\mu^1_{\La,\psi}(I_x)} 
\le \left(\frac{1-p_2}{1-p_1}\right)^{b/2}.
\]
Inequality \eqref{eqn:tech-lemma-start} follows as required.

\noindent(ii):
Inequality \eqref{eqn:tech-lemma-start} follows from
\eqref{rccomp}--\eqref{element2} on noting that $b\le \d$.

\noindent(iii), (iv):
Under either set of conditions, 
$\mu^1_{\La,\psi^x}\lest \mu^2_{\La,\psi^x}$, implying that
$\mu^1_{\La,\psi^x} (I_x)\ge \mu^2_{\La,\psi^x} (I_x)$.
Inequality \eqref{eqn:tech-lemma-start} follows on noting
that $b\le \d$.
\end{proof}

\begin{proof}[Proof of Theorem \ref{thm:monotone-limits}]
These inequalities follow in the same way as for the
\rc\ measure (see \cite{G-RC}, Section 4.3)
using the monotonicity of $\Phi^\l_{\La,a,p,q}$
for $\l=\bzero,\bone$.
\end{proof}

\begin{proof}[Proof of Theorem \ref{phimon}]
Let $C\subseteq \Psi\times\Om$ be an increasing cylinder event.
By the coupling of Section \ref{section:3},
$$
\phi^\l_{\La,a,p,q}(C) = \Phi^\l_{\La,a,p,q}(\mu^\l_{\La,\psi,p,q}(C_\psi)),
$$
where $C_\psi=\{\om\in\Om: (\psi,\om)\in C\}$ and $\mu_{\La,\psi,p,q}^\l$
is the \rc\ measure on $V(\psi)$ with boundary condition $\l$.
Now, $C_\psi$ is an increasing event in $\Om$, and therefore
$\mu^\l_{\La,\psi,p,q}(C_\psi))$ is increasing in $\psi$, $p$, and $\l$.
The claim follows by Theorem \ref{thm:comparison-inequality}(i).
\end{proof}

\begin{proof}[Proof of Theorem 
\ref{thm:finite-energy-Phi}]
Since $q\in[1,2]$, $\Phi^\l_{\La,a,p,q}$ is monotonic by Theorem
\ref{thm:positively-associated}. Since $J_x$ is increasing, a
lower bound for the conditional probability of $J_x$ is 
obtained by considering the situation in which
 all other vertices are closed.  In this case,
$x$ contributes $qa/(1-a)$ (\resp, 1) in 
\eqref{projrc2} when open 
(\resp, closed), and the lower bound follows.

An upper bound is obtained by considering the situation in which
$\l=\bone$, 
and all vertices other than $x$ are open and connected by open edges.
This time, $x$ contributes no more than
\[
r^\d q  \left(\frac{a}{1-a}\right)
\sum_{\om\in\{0,1\}^\d}
\,
\prod_{i=1}^\d 
\left(\frac{p}{1-p}\right)^{\om_i},
\]
when open, and $1$ when closed.
\end{proof}

\section{Stochastic orderings of edge-measures}
\label{EdgeComparisonSection}
Let $G=(V,E)$ be a finite graph, and let
$\phi_{a,p,q}$ be the \drc\ measure on the corresponding
sample space $\Psi\times\Om=\{0,1\}^V\times\{0,1\}^E$. 
Let $\EM_{a,p,q}$ denote the marginal measure of $\phi_{a,p,q}$
on the second component $\Om$,
\[
\EM_{a,p,q}(\om)=\sum_{\psi\in\Psi} \phi_{a,p,q}(\psi,\om),
\qq \om\in\Om.
\]
We first compare $\EM_{1,p_1,q_1}$ with $\EM_{a,p_2,q_2}$.

\begin{thm}
\label{EdgeComparisonThm}
Let $0< a_2 \le a_1 = 1$, $p_1,p_2\in (0,1)$, $q_1,q_2\in[1,\oo)$.
Let $r_i=\sqrt{1-p_i}$, and denote by 
$\EM_i$ the probability measure $\EM_{a_i,p_i,q_i}$.\\
\mbox{\rm(a)} If $q_2\le q_1$ and
\be
\label{eqn:EdgeComparisonCondition}
\frac{1-p_2}{p_2} (1 + 2 w_{\d}+ w_{\d} w_{\d-1}) \leq \frac{1-p_1}{p_1},
\ee
where $\d$ is the maximum vertex-degree of $G$ and
\[
w_j= \frac 1{q_2r_2^j}\left(\frac{1-a_2}{a_2}\right) ,
\qq j=0,1,2,\dots,\d,
\]
then $\EM_1\lest\EM_2$.
\\
\mbox{\rm(b)} If  $p_1\ge p_2$ and $q_1\le q_2$, then $\EM_1\gest\EM_2$.
\end{thm}

\begin{thm}
\label{EdgeComparisonThm2}
Let $0<a_1\le a_2<1$, $0<p_1\le p_2<1$, and $q\in[1,2]$. 
Then $\EM_{a_1,p_1,q} \lest \EM_{a_2,p_2,q}$.
\end{thm}

\begin{proof}[Proof of Theorem \ref{EdgeComparisonThm}]
(a) The quantity
\[
w_j(a,p,q)=\frac1{qr^j}\left(\frac{1-a}{a}\right) 
\]
may be viewed as follows. Let $(\psi,\om)\in\Theta$, and let $x\in
V$ be such that $\psi_x=0$.
Then
\be
\phi_{a,p,q}(\psi,\om)=\phi_{a,p,q}(\psi^x,\om) w_j,
\label{delclosed}
\ee
where $j=j(x,\psi)$ is the number of neighbours $u$ of $x$ such that
$\psi_u=1$.
Note that $w_j$ is increasing in $j$.

Suppose \eqref{eqn:EdgeComparisonCondition} holds.
We will show that the measures $\EM_i$ satisfy \eqref{Holleycond}. 
By the remark after \eqref{Holleycond}, it suffices to show that,
for $e,f\in E$ with $e\ne f$, and $\om\in\Om$ with $\om_e=0$, 
\begin{align}
\label{FKGEdgeComparison}
\EM_2(\om^{e,f})\EM_1(\om) &\ge \EM_1(\om^e)\EM_2(\om^f),\\
\label{FKGEdgeComparisonB}
\EM_2(\om^e)\EM_1(\om) &\ge \EM_1(\om^e)\EM_2(\om).
\end{align}
We will show \eqref{FKGEdgeComparison} only, the proof of \eqref{FKGEdgeComparisonB} 
is similar.  We may assume that $\om_f=0$.

Since $a_1=1$, $\EM_1$ is
the usual \rc\ measure on $G$ with parameters $p_1$ and $q$. Therefore,
\be\label{new10}
\EM_1(\om)=\EM_1(\om^e) \left(\frac{1-p_1}{p_1}\right) q_1^{k_1},
\ee 
where
\[
k_1=k(1,\om)-k(1,\om^e)=
\begin{cases}
1 & \text{if $e$ is an isthmus of the graph $(V,\eta(\om^e))$},\\
0 & \text{otherwise.}\\
\end{cases}
\]

For $\xi\in\Om$, let $K(\xi)=\{\psi: (\psi,\xi)\in\Th\}$ be the set
of compatible $\psi\in\Psi$.  
Let $e=\langle x,y\rangle$, and write
\[
B=\{\psi\in\Psi : \psi^{x,y}\in K(\om^{e,f}),\ \psi_x=\psi_y=0\}
\]
Then $K(\om^f)$ is the union of
\begin{romlist}
\item $\{\psi^{x,y}: \psi\in B\}$, and
\item $\{\psi^{x} : \psi\in B\}$ if $y$ is isolated in $\om^f$, and
\item $\{\psi^{y} : \psi\in B\}$ if $x$ is isolated in $\om^f$, and
\item $B$, if both $x$ and $y$ are isolated in $\om^f$.
\end{romlist}

Let $\psi\in B$. By \eqref{delclosed}, with $\phi_i=\phi_{a_i,p_i,q_i}$,
\begin{align*}
\phi_2(\psi^x,\om^f) 
&\begin{cases} = 0 &\text{if } (\psi^x,\om^f)\notin \Theta,\\
\le \phi_2(\psi^{x,y},\om^f)w_\d &\text{if } (\psi^x,\om^f)\in \Theta.
\end{cases}\\
&\le \phi_2(\psi^{x,y},\om^f)w_\d. 
\end{align*}
Similarly,
$$
\phi_2(\psi,\om^f)\le \phi_2(\psi^x,\om^f)w_\d
\le \phi_2(\psi^{x,y},\om^f)w_\d w_{\d-1}.
$$
Also,
$$
\phi_2(\psi^{x,y},\om^f)=\phi_2(\psi^{x,y},\om^{e,f})\left(\frac{1-p_2}{p_2}\right)
q_2^{k_2}, \qq \psi\in B,
$$
where
\[
k_2=k(\psi^{x,y},\om^f)-k(\psi^{x,y},\om^{e,f})
\le k_1.
\]
Therefore,
for $\psi\in B$,
\begin{multline*}
\phi_2(\psi^{x,y},\om^f)+\phi_2(\psi^x,\om^f)+\phi_2(\psi^y,\om^f)
+\phi_2(\psi,\om^f)\\
\le
\phi_2(\psi^{x,y},\om^{e,f}) \left(\frac{1-p_2}{p_2}\right) q_2^{k_2} ( 1 + 2w_\d
+ w_\d w_{\d-1}).
\end{multline*}
We sum over $\psi\in B$ and use \eqref{eqn:EdgeComparisonCondition}
and \eqref{new10} to find as required that
\begin{align*}
\EM_1(\om^e)\EM_2(\om^f) &\le \EM_1(\om^e)\EM_2(\om^{e,f}) 
\left(\frac{1-p_2}{p_2}\right)
q_2^{k_2} (1+2 w^\d + w_\d w_{\d-1})\\
&\le \EM_1(\om^e)\EM_2(\om^{e,f}) \left(\frac{1-p_1}{p_1}\right)q_1^{k_1}\\
&= \EM_1(\om)\EM_2(\om^{e,f}).
\end{align*}

\noindent(b)
The proof is similar but easier to that of (a), and is omitted.
\end{proof}

\begin{proof}[Proof of Theorem \ref{EdgeComparisonThm2}]
Write
$\phi_i=\phi_{a_i,p_i,q}$. For any increasing event $A\subseteq\Om$,
$$
\EM_1(A) = \phi_1(\Psi\times A)
=\phi_1(\phi_1(\Psi\times A\mid \psi))
=\Phi_1(\mu_{\psi,p_1,q}(A)),
$$
where $\mu_{\psi,p,q}$ denotes the \rc\ measure on $(V,E_\psi)$
with parameters $p$ and $q$. Now, $\mu_{\psi,p_1,q}\lest \mu_{\psi,p_2,q}$,
and $\mu_{\psi,p_2,q}(A)$ is non-decreasing in $\psi$. It follows
by Theorem \ref{thm:comparison-inequality}(i)
that $\Psi_1(A) \le \Psi_2(A)$ as required.
\end{proof}

\section{Infinite-volume measures}
\label{section:infiniteVolumeLimit}
There are two ways of moving to infinite-volume measures
on the lattice $\Ld=(\ZZ^d,\Ed)$,
namely by passing to weak limits, and by the 
Dobrushin--Lanford--Ruelle (DLR)
formalism. The associated theory is
standard for the \rc\ model, and the same arguments
are mostly valid for the \drc\ model. We shall not repeat
them here, but refer the reader to 
\cite{GrimmettUniquenessRandomClusterMeasures,G-RC}
for the details.

A subset of $\ZZ^d$ of the form $V_{a,b}=
\prod_{i=1}^d [a_i,b_i]$
is called a {\it box\/}, 
and the associated region is denoted by $\La_{a,b}$ and called
a {\it box-region\/}. Write $\sB$ for the set of all
box-regions of $\Ld$. For a sequence $\La_n$ of box-regions,
we write $\La_n\uparrow \LL^d$ if their vertex-sets increase
to $\ZZ^d$. Let $\Psi=\{0,1\}^{\Zd}$, $\Omega=\{0,1\}^{\Ed}$,
and let $\Theta$ be the set of all compatible pairs
$(\psi,\om)\in\Psi\times\Om$.

We begin with a consideration of vertex-measures.
Let $a,p\in(0,1)$ and $q\in(0,\oo)$,
and let $\sG$ denote the $\s$-field generated by the cylinder events of 
$\Psi=\{0,1\}^{\ZZ^d}$.
A probability measure $\Phi$ on $(\Psi,\sG)$ 
is called a {\it limit vertex-measure\/} with parameters $a$, $p$, $q$,
if, for some $\l\in\Th$, $\Phi$ is an 
accumulation point of the family 
$\{\Phi^{\l}_{\La,a,p,q}:\La\in\sB\}$. 
Let $\Wapq$ denote the set of all such measures, and 
$\coWapq$ its closed convex hull.
It is standard by compactness that $\Wapq$ is non-empty for all $a$, $p$, $q$.

We suppose henceforth that $q\in[1,2]$, so that we are within the
domains of validity of the comparison and
positive-correlation theorems of Sections 
\ref{BoundaryConditionsStochasticOrderings} and \ref{EdgeComparisonSection}.
Arguing as for \rc\ measures, any $\Phi\in\Wapq$ is positively
associated, and any $\Phi\in\coWapq$ has the finite-energy
property and satisfies the bounds of
Theorem \ref{thm:finite-energy-Phi}.

We may identify two special members of $\Wapq$ as follows. Let $\bzero
=(0,0)\in\Psi\times\Om$ and 
$\bone=(1,1)$. By positive-association in the usual way, the 
(monotonic) weak limits
$$
\Phi^b_{a,p,q} = \lim_{\La\uparrow\Ld} \Phi^b_{\La,a,p,q},\qq b=\bzero, \bone,
$$
exist. Furthermore, $\Phi^\bzero_{a,p,q}$ and $\Phi^\bone_{a,p,q}$
are automorphism-invariant (that is, invariant with respect to automorphisms
of $\LL^d$), and are extremal in that
\be
\Phi^\bzero_{a,p,q} \lest \Phi \lest \Phi^\bone_{a,p,q},\qq \Phi\in\coWapq.
\label{uniqm}
\ee
As in \cite{BenjamininiLyonsPeresSchrammUniformSpanningForests} (see also 
Section 4.3 of \cite{G-RC}), $\Phi^\bzero_{a,p,q}$
and $\Phi^\bone_{a,p,q}$ are tail-trivial, and are ergodic with respect
to the group $\ZZ^d$ of translations of $\Ld$.
Since  they have the finite-energy property, 
the number $I$ of infinite
open clusters satisfies either $I=0$ or $I=1$,
$\Phi^b_{a,p,q}$-a.s.\ ($b=\bzero,\bone$),
see \cite{MR990777,G-RC}. As noted after Theorem \ref{thm:monotone-limits},
the boundary conditions
$b=\bzero, \bone$ contain information concerning both vertex
and edge configuration off $\La$, but
only the external edge configuration is in fact relevant.

We shall perform comparisons in Sections \ref{ptransitions} 
and \ref{phasediagram}
involving these two extremal measures, 
and towards that end we note that, by weak convergence,
they satisfy the infinite-volume equivalents of
Theorem \ref{thm:comparison-inequality}.

The next two theorems concern the existence of the infinite-volume
limits for the \drc\ measure and the BCP measure, 
when $1\le q\le 2$. Here is a point of notation. Let
$\La=(V,E)$ be a box-region of $\Ld$, $q\in\{1,2\}$, and $s\in\{0,1,\dots,q\}$.
We write $\pi_{\La,K,\De,q}^s$ for the BCP measure on
$\La$ with boundary condition $s$. The boundary condition $s=0$ corresponds
to the free boundary condition. 
For $\psi\in\Psi$, $\mu_{\La,\psi,p,q}^b$ denotes
the \rc\ measure on $(V^+,E_\psi)$ with parameters $p$, $q$ and boundary condition $b$.
Similarly, $\mu_{\psi,p,q}^b$ denotes the corresponding \rc\ measure on
the infinite graph $(\Zd_\psi,\Ed_\psi)$. 
We write $\sH=\s(\sG\times\sF)$ for the product $\s$-field of $\Psi\times\Om$.
For $A\in\sH$ and $\psi\in\Psi$, let $A_\psi$ denote the section
$\{\om\in\Om:(\psi,\om)\in A\}$. For $B\subseteq \ZZ^d$, we write $B\lra \oo$ if there
exists $b\in B$ that is the endvertex of an infinite
open path of the lattice.

Let $\Vapq$ denote the
set of all weak-limit \drc\ measures with parameters $a$, $p$, $q$,
and let $\coVapq$ denote its closed convex hull. 
It is standard by compactness that
$\Vapq\ne\es$ for $a,p\in(0,1)$ and $q\in(0,\oo)$, and
by taking a Ces\`aro average of measures that $\coVapq$ 
contains some translation-invariant measure.
By part (a) of the next theorem, $\phi^b_{a,p,q}\in \Vapq$ when
$q\in[1,2]$. 

\begin{thm} \label{infvolthm1}
Let $a,p\in(0,1)$, $q\in[1,2]$, and $b\in\{\bzero,\bone\}$.\\
\mbox{\rm(a)} The limit
\drc\ measure $\phi_{a,p,q}^b = \lim_{\La\uparrow\Ld}\phi_{\La,a,p,q}^b$
exists and satisfies
$$
\phi_{a,p,q}^\bzero(A) = \Phi^\bzero_{a,p,q}(\mu_{\psi,p,q}^0(A_\psi)),
\qq A\in \sH,
$$
with a similar equation for the boundary condition $\bone$.\\
\mbox{\rm(b)} The $\phi^b_{a,p,q}$ are stochastically
increasing in $a$ and $p$, and $\phi^\bzero_{a,p,q}
\lest \phi \lest \phi^\bone_{a,p,q}$ for $\phi\in\coVapq$.\\
\mbox{\rm(c)} We have that
$$
\phi_{\La,a,p,q}^\bone(0\lra\oo) \to \phi_{a,p,q}^\bone(0\lra\oo)
\qq\text{as } \La\uparrow\Ld.
$$
\mbox{\rm(d)} The number $L(\om)$ of infinite open clusters
of $\om\in\Om$ satisfies: either $\phi^b_{a,p,q}(L=0)=1$
or $\phi^b_{a,p,q}(L=1)=1$.
\end{thm}

\begin{thm}\label{infvolthm2}
Let $K\in[0,\oo)$, $\De\in\RR$, and $q\in \{1,2\}$.
The limit BCP measure 
$\pi_{K,\De,q}^s=\lim_{\La\uparrow\Ld}\pi_{\La,K,\De,q}^s$
exists, for $s=0,1,\dots,q$.
\end{thm}

The proofs are deferred to the end of this section.
We recall from Section 3 the `usual' coupling of the \drc\ and BCP measures,
and we shall see in the proof of the last theorem
that the equivalent coupling is valid for
the infinite-volume measures.

The limit measures $\phi_{a,p,q}^b$ are
automorphism-invariant and have the finite-energy property,
the proofs follow standard lines and are omitted.
Similarly, the $\phi_{a,p,q}^b$ satisfy
the comparison inequalities of
Theorems \ref{EdgeComparisonThm} and \ref{EdgeComparisonThm2}.

We shall consider
also the set of DLR measures. 
Let $\sT_\La$ be the sub-$\s$-field of $\sH$ generated by the states
of vertices and edges not belonging to the region $\La$.
A probability measure on $(\Psi\times\Om,\sH)$
is called a {\em \drc\ measure} with parameters $a$, $p$, $q$ 
if, for every $A\in\sH$ and every region $\La$,
\[
\phi(A \mid \sT_\La) (\t) = \phi^\t_{\La,a,p,q} (A) \qq 
\text{for }\phi\text{-a.e. }\t\in\Psi\times\Om.
\]
The set of such measures is denoted by $\Rapq$. 
One way of showing that $\Rapq\ne \es$ is to prove that
some measure in $\coVapq$ belongs to $\Rapq$. The following
theorem may be proved exactly as for \rc\ measures, 
see \cite{GrimmettUniquenessRandomClusterMeasures,G-RC}.

\begin{thm}
\begin{romlist}
\item Let $a,p\in(0,1)$ and $q\in(0,\oo)$.
If $\phi\in\coVapq$ and $\phi$  is such that $\phi(L\in\{0,1\})=1$,
then $\phi\in\Rapq$.
\item $\Rapq\ne\es$ for $a,p\in(0,1)$, $q\in(0,\oo)$.
\item
Let $a,p\in(0,1)$ and $q\in[1,2]$. Then $\phi^b_{a,p,q}
\in\Rapq$ for $b=\bzero,\bone$.
\end{romlist}
\end{thm}

Finally, we indicate how the convexity of the partition function
may be used to show the uniqueness of certain infinite-volume measures.
The proof follows \cite{GrimmettUniquenessRandomClusterMeasures}, 
which in turn used the method of \cite{LM-L}.

\begin{thm}\label{uniqevertexm} Let $q\in[1,2]$. 
\\
\mbox{\rm(a)} For $p\in(0,1)$, the set of points $a\in(0,1)$ at which
$|\Wapq|\ge 2$ is countable.
\\
\mbox{\rm(b)} If $q\in\{1,2\}$, the set of pairs $(a,p) \in (0,1)^2$ at which
$|\Vapq|\ge 2$ may be covered by a countable family of rectifiable curves of $\RR^2$.
\end{thm}

\begin{proof}[Proof of Theorem \ref{infvolthm1}]
(a) For simplicity in the following proofs, 
we shall suppress reference to the parameters.
Consider first the boundary condition $\bzero$.
Let $A\subseteq \Om$ and $B\subseteq \Psi$ be increasing cylinder events, and
let $U\subseteq\ZZ^d$ be a finite set such that $A$ and $B$
are defined in terms of
the states of vertices in $U$ and of edges joining members of $U$. 
By the discussion in Section \ref{BoundaryConditionsStochasticOrderings}, 
\be\label{new50}
\phi^\bzero_{\La}(A\times B) =
\Phi^\bzero_{\La}(1_A(\psi)\mu_{\La,\psi}^0(B)).
\ee
Since $\sH$ is generated by the set of such events $A \times B$,
it suffices to show that
\be\label{new51}
\lim_{\La\uparrow\Ld} \phi_{\La}^\bzero(A\times B) = \Phi^\bzero(1_A(\psi)\mu_{\psi}^0(B)).
\ee
Let $\La'=(V',E')$, $\La''$ be box-regions such that $\La'\subseteq \La\subseteq \La''$
and $U\subseteq V'$.
By \eqref{new50} and the
monotonicity of  $\Phi_\La^\bzero$ in $\La$, and of $\mu_{\La,\psi}^0$ in $\La$
and $\psi$,
$$
\Phi^\bzero_{\La}(1_A(\psi)\mu_{\La',\psi}^0(B)) \le
\Phi^\bzero_{\La}(A\times B) \le
\Phi^\bzero_{\La''}(1_A(\psi)\mu_{\La,\psi}^0(B)).
$$
Take the limits as $\La'',\La, \La'\uparrow\Zd$ in that order,
and use the bounded convergence theorem to obtain 
\eqref{new51}.
A similar argument holds with boundary condition 
$\bone$, and with the inequalities reversed.

\noindent
(b) The necessary properties of monotonicity follow by Theorem \ref{phimon}.

\noindent
(c) This follows the proof of the corresponding statement for
\rc\ measures, see \cite{MR939480,G-RC}, using part (a) and the representation
\eqref{new50} with boundary condition $\bzero$ replaced by $\bone$.

\noindent
(d) The proof relies on the automorphism-invariance and the finite-energy property
of the marginal measure of $\phi_{a,p,q}^b$ on $\Om$. This
follows standard lines and is omitted.
\end{proof}

\begin{proof}[Proof of Theorem \ref{infvolthm2}]
Consider first the case $s=0$. Let $\La^-$ be the graph
obtained from the box-region $\La=(V,E)$ by 
removing those edges that do not have both endvertices in $V$.
Let $\mu$ be the coupled measure of Theorem \ref{thm:couple-BC-DRC} for 
$\La^-$, having marginal measures $\pi^0_\La=\pi^0_{\La,K,\De,q}$ 
and $\phi^\bzero_\La = \phi^\bzero_{\La,a,p,q}$, where $a$, $p$ satisfy \eqref{apdef}.

Let $U\subset \ZZ^d$ be finite, $\tau\in\Si=\{0,1,2,\dots,q\}^{\Zd}$, and 
let $\Si_{U,\tau}$ be the BCP cylinder event $\{\s\in\Si:\s_u=\tau_u
\text{ for } u\in  U\}$. Let $A=A_{U,\tau}$ be the set
of $\t=(\psi,\om)\in\Theta$ that are compatible with $\Si_{U,\tau}$, that is,
$A$ is the set of $\t$ such that:
\begin{romlist}
\item $\forall u \in U$, $\tau_u=0$ if and only if $\psi_u=0$, and
\item $\forall u,v \in U$, $\tau_u \ne \tau_v$ only if $u$ and $v$ are not 
$\om$-connected in $\Ld$.
\end{romlist}
For given $\t\in A$, let $l(\t)$ be the number of open clusters that 
intersect $U$. By the second observation after Theorem 
\ref{thm:couple-BC-DRC}, subject to a slight abuse of notation, 
if $V\supseteq U$,
\be\label{phipi0}
\pi^0_{\La} (\Si_{U,\tau}) = 
\phi^\bzero_{\La}\bigl(1_{A}(\t)q^{-l(\t)}\bigr).
\ee
Now, $\phi^\bzero_{\La}\Rightarrow \phi^\bzero$
as $\La\uparrow\Ld$ and, by Theorem \ref{infvolthm1}(d), 
the random variable $1_{A}(\t)q^{-l(\t)}$
is $\phi^\bzero$-a.s.\ continuous. Therefore,
$$
\lim_{\La\uparrow\Ld} \pi^0_{\La} (\Si_{U,\tau}) =
\phi^\bzero\bigl(1_{A}(\t)q^{-l(\t)}\bigr).
$$

Suppose now that $s\in\{1,2,\dots,q\}$.
Let $\mu$ be the coupled measure of Theorem \ref{thm:couple-BC-DRC} 
on the graph $(V^+,E)$, and let $\mu^s$ denote
the measure $\mu$ conditioned on the 
event that $\s_x=s$ for all $x\in\pd\La$.

The marginal of $\mu^s$ on $\Si_V=\{1,2,\dots,q\}^V$ 
is the measure $\pi^s_\La = \pi^s_{\La,K,\De,q}$,
the marginal on $\Psi_V\times\Om_E=\{0,1\}^V\times\{0,1\}^E$ 
is $\phi^\bone_\La = \phi_{\La,a,p,q}^\bone$.
The conditional measure of $\mu^s$ on $\Si_V$, given the
pair $(\psi,\om)\in \Psi_V\times\Om_E$, is that obtained as follows:
\begin{letlist}
\item $\forall v\in V$, the spin at $v$ is $0$ if and only if $\psi_v=0$,
\item the spins are constant on each given open cluster,
\item the spins on any open cluster intersecting $\pd\La$ are equal to $s$,
\item the spins on the other open clusters are independent 
and uniformly distributed
on the set $\{1,2,\dots,q\}$.
\end{letlist}
Equation \eqref{phipi0} becomes
\be\label{phipis}
\pi^s_{\La} (\Si_{U,\tau}) =
\phi^\bone_{\La}\bigl(1_{A}(\t)q^{-f(\t)}\bigr),
\ee
where $f(\t)$
is the number of finite open clusters that intersect $U$. 
[Recall that $\phi^\bone_\La$ has support $\Theta_\La^\bone$.]
We may write  $f(\t) = l(\t) - N(\om)$ where
$N=N(\om)$ is the number of infinite open clusters of $\om\in\Om$
that intersect $U$.
Clearly, $N\in\{0,1\}$ for $\t=(\psi,\om)\in\Th_\La^\bone$, so that
\begin{align}
\pi^s_{\La} (\Si_{U,\tau}) &=
\phi_\La^\bone(1_Aq^{N-l})\nonumber\\
&= \phi_\La^\bone(1_Aq^{-l}) + (q-1)\phi_\La^1(1_A1_{\{U\lra\pd\La\}}
q^{-l}).
\label{new75}
\end{align}
Now, $1_Aq^{-l}$ is $\phi^\bone$-a.s.\ continuous by
Theorem \ref{infvolthm1}(d), so that
\be\label{new53}
\phi_\La^\bone(1_Aq^{-l}) \to \phi^\bone(1_Aq^{-l})\qq\text{as } \La\uparrow\Ld.
\ee
It may be proved in a manner very similar to the proof
of Theorem \ref{infvolthm1}(c) that
\be\label{new76}
\phi_\La^\bone(1_A1_{\{U\lra\pd\La\}}q^{-l})
\to \phi^\bone(1_A1_{\{U\lra\oo\}}q^{-l})
\qq\text{as } \La\uparrow \LL^d.
\ee
By \eqref{new75}--\eqref{new76} and Theorem \ref{infvolthm1}(d),
$$
\pi^s_{\La} (\Si_{U,\tau}) \to
\phi^\bone(1_Aq^{-f})
\qq\text{as } \La\uparrow\LL^d,
$$
and the proof is complete.
\end{proof}

\begin{proof}[Proof of Theorem \ref{uniqevertexm}]
(a) Let $\La=(V,E)$ be a region in $\Ld$ with graph $\La^+=(V^+,E)$. 
Let $a,p\in(0,1)$ and $q\in[1,\oo)$.
Consider the normalizing constant 
$Z^\l_\La=\ZDRC_{\La,\l,a,p,q}$ of the \drc\ measure
on $\La$ with boundary condition $\l$.
Let the vectors $(a,p)$ and $(K,\De)$ satisfy \eqref{apdef}. 
By \eqref{new34}, we may write
$$
Z^\l_\La=\sum_{\theta=(\psi,\om)\in\Theta_\La^\l} 
r^{|E_\psi|} q^{k(\t,\La)} 
e^{-\De|V_\psi|}\left(\frac p{1-p}\right)^{|\eta(\om)\cap E|}.
$$
By a standard argument using subadditivity in $\La$, 
see \cite{GrimmettUniquenessRandomClusterMeasures,G-RC},
the limit
$$
G(\De,p, q) = \lim_{\La\uparrow\Ld} \left\{\frac 1{|V|} \log Z_\La^\l\right\}
$$
exists and is independent of $\l$. 
The function $G$ is termed {\it pressure\/}.

It is easily seen that
\begin{align}
\label{diff1}
\frac\pd{\pd \De} \log Z_\La^\l &= -\phi^\l_{\La}(|V_\psi|),\\
\label{diff2}
\frac{\pd^2}{\pd^2 \De} \log Z_\La^\l &= \Var(|V_\psi|),
\end{align}
where $\Var$ denotes variance with respect to $\phi^\l_{\La,a,p,q}$.
Since variances are non-negative, $G(\De,p,q)$ is a convex
function of $\De$.
Hence, for fixed $p$, $q$, the set of points $\De$ of
non-differentiability of $G$ is countable (that is, either finite
or countably infinite).
Wherever $G$ is differentiable, its derivative is the limit
as $\La\uparrow\Ld$ of the 
derivative of $|V|^{-1}\log Z^\l_\La$. This implies
in turn that
\[
\lim_{\La\uparrow\Ld} \frac{1}{|V|} \phi^\bzero_{\La}(|V_\psi|)
=\lim_{\La\uparrow\Ld} \frac{1}{|V|} \phi^\bone_{\La}(|V_\psi|),
\]
so that
$\Phi^\bzero(J_x)=\Phi^\bone(J_x)$ for $x\in \ZZ^d$, where
$J_x$ is the event that $x$ is open. 
The claim follows by \eqref{uniqm}
and a standard `FKG' coupling (see, for example, Prop.\ 4.6 of \cite{G-RC}).

\noindent
(b) When $q\in\{0,1\}$, we work with the constant
$\ZBCP=\ZBCP_{\La,K,\De,q}$ of \eqref{BCPdef}. By the form of \eqref{BCPdef},
$\ZBCP_{\La,K,\De,q}$ is a convex function of the pair $(K,\De)$.
By \eqref{eqpartfns} and the coupling of Chapter 3,
\begin{align}
\label{diff1b}
\frac\pd{\pd K} \log Z_{\La}^b &=
 \pi_{\La}^s \biggl(-|E_\s|+ 2 \sum_{e\in E}
\d_e(\s)\biggr)\nonumber\\
&=\phi_{\La}^b \biggl(-|E_\psi|+ \frac2p \sum_{e\in E} \om(e)\biggr),
\end{align}
where $s=s(b)$ satisfies $s(\bzero)=0$, $s(\bone)=1$.
By Theorem 8.18 of
\cite{Fal} or Theorem 2.2.4 of \cite{Schn},
the set of points of $(0,1)^2$ at which
$G$ is not differentiable (when viewed as function of
$(a,p)$) may be covered by a countable
collection of rectifiable curves.
Suppose $G$ is differentiable at the point $(a,p)$. By part (a),
$\phi^\bzero(J_x) = \phi^\bone(J_x)$ for
$x\in \ZZ^d$ and, in particular,
$|E_\psi|/|V|$ has the same (almost-sure and $L^1$) limit
as $\La\uparrow\LL^d$
under either boundary condition.
Therefore, by \eqref{diff1b},
\[
\lim_{\La\uparrow\Ld} \frac{1}{|V|} \phi^\bzero_{\La}(|\eta(\om)\cap
E|)=\lim_{\La\uparrow\Ld} \frac{1}{|V|} \phi^\bone_{\La}(|\eta(\om)\cap E|),
\]
so that, by translation invariance, $\phi^\bzero(J_e)=\phi^\bone(J_e)$
for $e\in\Ed$, where $J_e$ is the event that $e$ is open. 
The claim now follows by Theorem \ref{infvolthm1}(b),
as in part (a).
\end{proof}

\section{Phase transitions}\label{ptransitions}
Let $d\ge 2$, $q\in[1,2]$, and consider the `wired' \drc\ measure
$\papq$ on $\Ld$. Several transitions occur as $(a,p)$ increases
from $(0,0)$ to $(1,1)$, and each gives rise to a 
`critical surface' defined as follows. 

Let $\Pi$ be a monotonic
property of pairs $(\psi,\om)\in\Theta$ such that $\papq(\Pi)\in\{0,1\}$
for all $a$, $p$. Let 
$$
R(\Pi)=\{(a,p)\in (0,1)^2: \papq(\Pi)=1\},\q
S(\Pi)=\overline{R(\Pi)}\cap \overline{R(\neg \Pi)},
$$
where
$\neg\Pi$ denotes the negation of $\Pi$. By Theorem \ref{infvolthm1}(a),
each $R(\Pi)$ is a monotonic subset of $(0,1)^2$ with respect to
the ordering $(a,p) \le (a',p')$ if
$a\le a'$ and $p\le p'$. The set 
$S(\Pi)$ ($=S(\neg\Pi)$) is called the `critical surface' for $\Pi$.

Of principal interest here are the following three properties:
\begin{romlist}
\item 
$\Picvc$, the property that there exists an infinite closed vertex-cluster,
\item
$\Piovc$, the property that there exists an infinite open vertex-cluster,
\item
$\Piec$, the property that there exists an infinite open edge-cluster.
\end{romlist}
It is easily checked that $\neg\Picvc$, $\Piovc$, $\Piec$ are increasing and 
satisfy the zero/one claim above. Furthermore, $\Piec \subseteq \Piovc$.

We do not know a great deal about the critical surfaces of the three 
properties above. Just as for percolation,
it can occur that $R(\Picvc) \cap R(\Piovc)\ne\es$ on any lattice
whose critical site-percolation-probability $\pcsite$ satisfies
$\pcsite<\frac12$, see remark (a) following \eqref{projrc2}.
When $d=2$ however, $R(\Picvc) \cap R(\Piovc)=\es$ by the
main theorem of \cite{MR942759}.

When $q=2$, the critical surfaces of these 
three properties mark phase transitions
for the \BC\ model. Consider the \BC\ measure $\pi_{K,\De,2}^1$ on $\Ld$, 
and let $a$, $p$ satisfy \eqref{apdef}. Then:
\begin{romlist}
\item
$R(\Picvc)$ corresponds to the existence of an infinite vertex-cluster of
spin $0$,
\item
$R(\Piovc)$ corresponds to the existence of an infinite vertex-cluster whose
vertices have non-zero (and perhaps non-equal) spins,
\item
$R(\Piec)$ corresponds to the existence of long-range order.
\end{romlist}
Statements (i)--(ii) are clear. Statement (iii) follows by
Theorems \ref{corrconn} and \ref{infvolthm1}(c), and the remark following Theorem \ref{infvolthm2}, 
on noting by \eqref{twoptcf} that
\be
\label{infcorrconn}
\pi^1_{K,\De,2}(\s_0=1) - \tfrac 12\pi_{K,\De,2}^1(\s_0\ne 0) =
\tfrac12\paptwo(0\lra \oo ).
\ee

Some numerical information may be obtained about the critical surfaces by use
of the comparison inequalities proved earlier in this paper. This
is illustrated in the next section, where we concentrate on the 
two-dimensional \BC\ model.

\begin{figure}[ht]
\begin{center}
\input{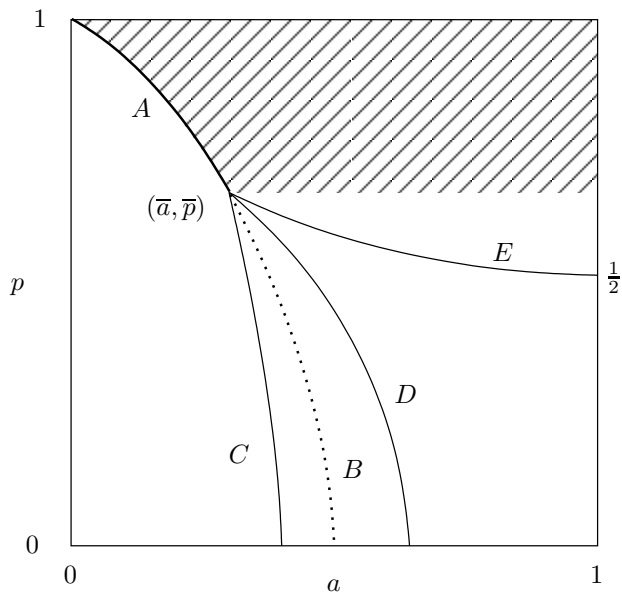}
\caption{The
phase diagram of the $q=1$ BCP model on the square lattice. 
The model may be transformed
into the Ising model with edge-interaction $J$ and external
field $h$, see \eqref{new36}. The arc $A\cup B$ with equation $a/(1-a)=1-p$
corresponds to $h=0$, and the points to its right (\resp, left)
correspond to $h>0$ (\resp, $h<0$).  There is a `tri-critical point' 
at $(\ol a,\ol p)$, see \eqref{olaolp2}, and the arc $A$ joining
this point to $(0,1)$
marks a line of first-order phase transitions. The region to
the left of $A\cup C$ is $R(\Picvc)$, and that to the right
of $A\cup D$ is $R(\Piovc)$. The hatched region lies in $R(\Piec)$,
and the lower boundary of $R(\Piec)$ is presumably as marked by $E$.}
\label{newfig3}
\end{center}
\end{figure}

This section closes with some notes on the BCP model on $\LL^2$ with $q=1$,
for use in Section \ref{phasediagram}.
As remarked in Sections \ref{section:3} and \ref{latticeDRC}, this 
model may be transformed into the Ising model
with edge-interaction $J=\frac14 K$ and external
field $h=K-\frac12 \De$, see \eqref{new36}. The phase diagram is
therefore
well understood and is illustrated in Figure \ref{newfig3}
with the parametrization $(a,p)$ of \eqref{apdef}.

Some remarks concerning Figure \ref{newfig3} follow. The existence
of the arcs $A$, $C$, $D$ follow by the established theory of the Ising model
with edge-interaction $J$ and external field $h$, see \cite{CNPR,Hig1,Hig2,Ru79}
for the case $h=0$. 
The arc $A$ corresponds
to $h=0$, $J>\Jc$, where $\Jc$ is the critical point of the zero-field model.
Consider the corresponding \rc\ model $\RC_{p}$ with edge-parameter $\pi=1-(1-p)^{-4}$
and cluster-weighting factor $2$. Then $\RC_p$ has (almost surely) an
infinite open cluster $I_p$ when $(a,p)\in A$.
As one deviates rightwards from $A$ with $p$ held constant (that is, in the direction
of positive $h$), 
the positive magnetic field attracts the
vertices in $I_p$, together with at least one half of the finite clusters
of $\RC_p$. Write $\P_{p,h}$ for the resulting set of $+1$ spins. 
By the previous remark, and recalling the conditional law of the 
zero-field Ising model given the \rc\ configuration, we deduce that
the bond percolation model on $\P_{p,h}$ with density $\pi$ ($<p$)
possesses an infinite edge-cluster. It follows that the hatched region of
Figure \ref{newfig3} lies in $R(\Piec)$.

Similarly, as one deviates leftwards from $A$ with $p$ held constant,
the resulting negative magnetic field attracts the vertices in $I_{p}$,
and an infinite closed vertex-cluster forms.

\section{The \BC\ phase diagram}\label{phasediagram}
\begin{figure}[ht]
\begin{center}
\includegraphics[angle=0,width=10cm]{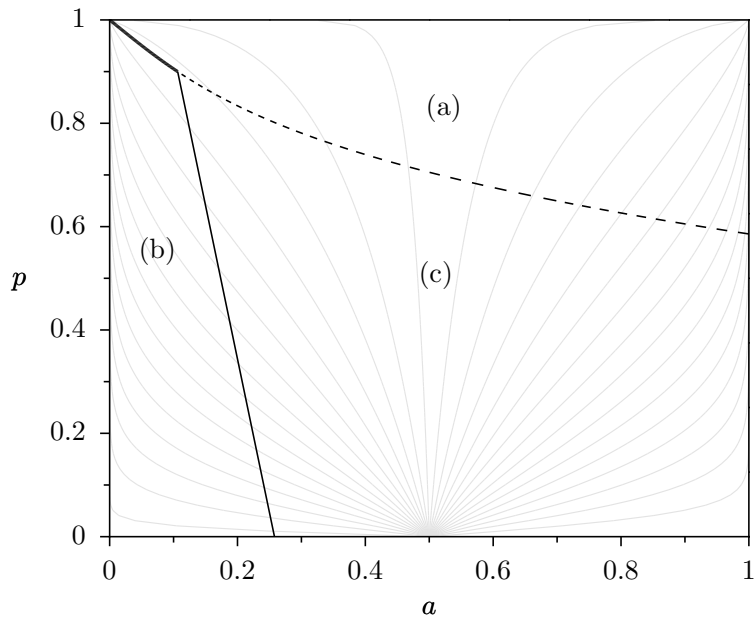}
\caption{The \BC\ phase diagram in two dimensions as proposed
by Capel. Note the three phases labelled (a), (b), (c) as in the text.
The boundary between (a) and (b) is thought to be a line of first-order
phase transitions, whereas that between (a) and (c) is expected
be a line of second-order transitions. The point at which the
the phases are expected to meet is termed the {\em tri-critical point}.
Moral support for such a phase diagram is provided by
the rigorously known $q=1$ diagram of Figure \ref{newfig3}.}
\label{fig1}\end{center}
\end{figure}

Throughout this final section, we consider the \BC\ model
on the square lattice $\LL^2$, and the associated \drc\ measure.
[Related but partial conclusions are valid similarly on $\LL^d$ with $d\ge 3$.]
The respective parameters are $K\in[0,\oo)$, $\De\in\RR$,
and the values $a$, $p$ given at \eqref{apdef}. The three putative phases
of the models are illustrated in Figure \ref{fig1}. 
We recall from the last section
the fact that, since $d=2$, $R(\Piovc)\cap R(\Picvc)=\es$.

The three regions of Figure \ref{fig1} are characterized as follows.
\begin{letlist}
\item The top region is $R(\Piec)$, in which the \drc\ measure possesses
(almost surely) an 
infinite open edge-cluster, and the \BC\ model 
has long-range order.
\item The left region is $R(\Picvc)$, in which
the measures possess 
an infinite vertex-cluster of zero states.
\item The central region is $R(\neg\Picvc)\cap R(\neg\Piec)$, in which
either all closed and open vertex-clusters are finite, or
there exists an infinite open vertex-cluster which is insufficiently small
to support an infinite
open edge-cluster. There is no long-range order.
\end{letlist}
In the more normal parametrization \eqref{BCHam} of the \BC\ model, there is
a parameter $\b$ denoting inverse-temperature, and one takes
$K=\b J$, $\De=\b D$. If we hold the ratio $D/J$ fixed and let $\b$
vary, the arc of corresponding pairs $(a,p)$ satisfies
\[
\frac{a}{1-a}=(1-p)^{D/2J}.
\]
As the ratio $D/J$ varies, such arcs are plotted in the gray lines
of Figure \ref{fig1}.

The region labelled (c) may be split into two sub-regions
depending on whether or not there exists an infinite open vertex-cluster.
We shall not pursue this distinction here.

A key prediction of Capel for this model is the existence
of a so-called tri-critical point where the 
three phases meet. 
The common boundary between the regions $R(\Picvc)$ and $R(\Piovc)$
is thought to be a line of first-order phase transitions. Based
on a mean-field analysis, Capel has
made the numerical proposals that the tri-critical point
lies on the line $a/(1-a)=(1-p)^{\frac23 \log 4}$, and that the line of
first-order transitions
arrives at the corner $(0,1)$ with the same gradient as the line
${a}/({1-a})=1-p$. 
The remaining boundary of $R(\Piec)$ is thought to mark a line of 
second-order phase transitions, and to meet the line $a=1$ at the point
$p=\sqrt 2/(1+\sqrt 2)$.

\begin{figure}[ht]
\begin{center}
\includegraphics[angle=0,width=10cm]{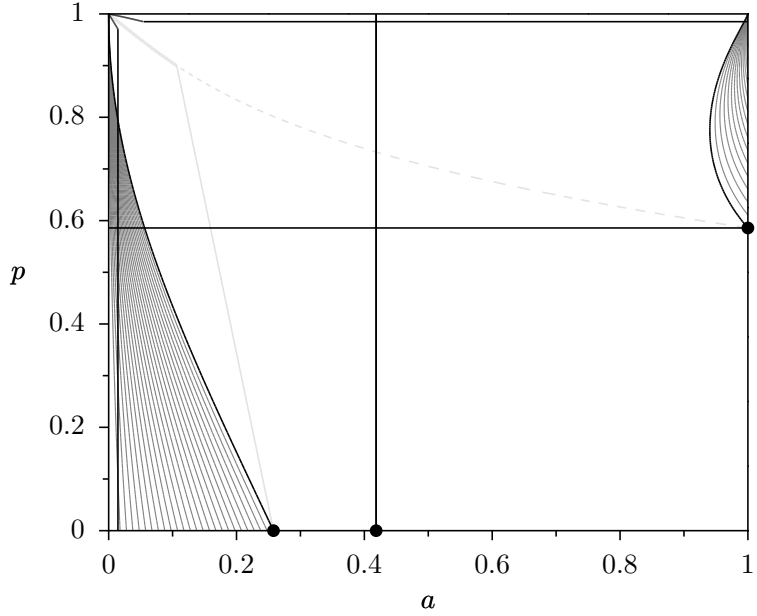}
\caption{Regions of the phase diagram of the Blume--Capel
model on $\LL^2$ about which one may make
rigorous statements on the basis of comparisons with other models.
The three points referred to in (i)--(iii) are marked. The narrow vertical
strip along the $p$-axis is a subset of $R(\Picvc)$, and the horizontal strip
along the line $p=1$ is a subset of $R(\Piec)$; 
see the comments around \eqref{new64} and \eqref{new65}.}
\label{fig2}
\end{center}
\end{figure}

The $q=2$ \rc\ (Ising) measure on $\LL^2$ has 
critical point $p={\sqrt{2}}/(1+\sqrt{2})$,
which for numerical clarity we shall approximate by $0.586$. 
Site percolation on $\LL^2$ has critical probability
$\pcsite$, to which we shall approximate with the value $0.593$.
Figure \ref{fig2} indicates certain regions of the phase diagram 
about which we may make precise observations.

For three special vectors $(a,p,q)$, the corresponding
\drc\ measure $\phi^\bone_{a,p,q}$ provides information concerning the
phase diagram. These vectors are given as follows.
For simplicity, we shall refer to the comparison theorems for
measures on finite graphs; the corresponding inequalities for infinite-volume
measures are easily seen to hold, see Section \ref{section:infiniteVolumeLimit}.
\begin{romlist}
\item
{\em The triple $a=1$, $p={\sqrt{2}}/({1+\sqrt{2}})\approx 0.586$, $q=2$.}
The corresponding $\phi^\bone_{1,p,2}$ is a critical \rc\ measure.
By Theorem \ref{EdgeComparisonThm}(a), the shaded region to the right of the 
given curve 
joining $(1,0.586)$ to $(1,1)$ lies within $\Piec$. 
The corresponding \BC\ models have long-range order.
By Theorem \ref{EdgeComparisonThm}(b), no point below
the horizontal line $p=0.586$ lies in $\Piec$,
and the corresponding \BC\ models do not have long-range order.
\item
{\em The triple $a=({1-\pcsite})/({1+\pcsite})\approx 0.26$, $p=0$, $q=2$.}
To the left of this point on 
the horizontal line $p=0$, the vertex-measure $\Phi^\bone_{a',0,2}$ 
is a product measure with density $2a'/(1+a')$ and 
possessing (almost surely) an infinite cluster of closed vertices.
By Theorem \ref{thm:comparison-inequality}(ii), $\Phi^\bone_{a',0,2}$  
dominates the vertex-measures on the given arc joining $(a',0)$ to $(0,1)$. 
The interior of the shaded area is thus a subset of $R(\Picvc)$,
and the corresponding BCP measures possess (almost surely) an 
infinite cluster of $0$-spin vertices.
\item
{\em The vector $a=\pcsite/(2-\pcsite)\approx 0.42$, $p=0$, $q=2$.}
To the right of this point on the horizontal line $p=0$, 
the vertex-measure $\Phi^\bone_{a',0,2}$ is a supercritical product 
measure with an infinite open vertex-cluster.
It follows by Theorem \ref{thm:comparison-inequality}(i) that
the interior of the region to the right of the vertical line 
$a=\pcsite/(2-\pcsite)$
lies in $R(\Piovc)$.
\end{romlist}

Finally, we shall make comparisons involving the \drc\ model with parameters
$(a,p,2)$ and the $q=1$ models lying on the arc $A$ of Figure \ref{newfig3}.
Let $\ol a\approx 0.029$ and 
$\ol p \approx 0.971$ be given by \eqref{olaolp2}, and consider
the BCP model with parameters $(a_2,p_2,1)$ where
$a_2/(1-a_2)=1-p_2$. 
Take $q_2=1$ and
$q_1=2$ in Theorem \ref{thm:comparison-inequality}(iii) to find
that: if $(a,p)\in(0,1)^2$ satisfies
$$
2\left(\frac a{1-a}\right)(1-p)^2 < (1-p_2)^3
$$
for some $p_2\ge \max\{p,\ol p\}$, then $(a,p)\in R(\Picvc)$.
This holds in particular if 
\be\label{new64}
\frac{2a}{1-a} < 1-p\q\text{and}
\q p>\ol p.
\ee
Taken in conjunction with Theorem \ref{thm:comparison-inequality}(i), this
implies that the narrow vertical strip marked along the $p$-axis of Figure \ref{fig2} 
is a subset of $R(\Picvc)$.

Secondly, take $q_1=1$, $q_2=2$  in
Theorem \ref{thm:comparison-inequality}(iv) to find similarly that:
if $(a,p)\in(0,1)^2$ satisfies
$$
2\left(\frac a{1-a}\right)(1-p)^2 > \left(1-\frac p{2-p}\right)^3
\q\text{and}\q \frac p{2-p} > \ol p,
$$
then $(a,p)\in R(\Piovc)$. This occurs if 
\be\label{new65}
\frac {2a}{1-a} > \frac{8(1-p)}{(2-p)^3}\q\text{and}
\q p > \frac{2\ol p}{1+\ol p}.
\ee

We indicate next that $(a,p)\in R(\Piec)$ whenever
\eqref{new65} holds. Assume \eqref{new65}. By
Theorem \ref{thm:comparison-inequality}(iv),
$\Phi^\bone_{a,p,2} \gest \Phi^\bone_{a_1,p_1,1}$ where 
$a_1/(1-a_1)=1-p_1$ and $p_1=p/(2-p)>\ol p$.
Since the inequalities of \eqref{new65} are strict, we may replace
$a_1$ by $a_1+\epsilon$ for some small $\epsilon>0$, and we deduce that 
$\Phi^\bone_{a,p,2}$ dominates (stochastically) the law,
$\mu_J^+$ say,  of the set $S$ of
$+$-spins of the infinite-volume Ising model with zero external field, edge-interaction
$J=-\frac 18\log(1-p)>\Jc$, and $+$ boundary condition. 
Recalling the coupling between the Ising model and the \rc\ model,
the critical probability $\pcb(S)$ of bond percolation on
$S$ satisfies $\pcb(S) < \pi$, $\mu_J^+$-a.s., where $\pi$
is the `effective' edge-parameter of the \rc\ model $\RC_{p_1}$ given by
$$
(1-\pi)^4 = 1-p_1 = 1-\frac p{2-p}. 
$$
The \rc\ measure
with parameters $p$, $2$ on the graph induced by the open vertex-set of $\LL^2$
dominates (stochastically) the product measure with intensity
$p_1=p/(2-p)$. Since $p_1\ge \pi$, there exists an infinite
open edge-cluster, $\phi^\bone_{a,p,2}$-a.s.
That is, $(a,p)\in R(\Piec)$ if \eqref{new65} holds.
This implies as above that the narrow horizontal strip marked along the
line $p=1$ in Figure \ref{fig2} is a subset of $R(\Piec)$.

\section*{Acknowledgements}
We thank Aernout van Enter for his advice on the literature.
The first author acknowledges financial support from the 
Engineering and Physical Sciences Research Council 
under a Doctoral Training Award to the University of Cambridge.

\bibliography{blume}
\bibliographystyle{plain}

\end{document}